\documentclass[reqno,11pt]{amsart}

\usepackage[utf8]{inputenc}
\usepackage[T1]{fontenc}
\usepackage{lmodern}

\usepackage[numbers,sort&compress]{natbib} 
\usepackage{hyperref}
\usepackage{enumitem}
\usepackage{graphicx}
\usepackage[labelfont={rm}, justification=centering]{subfig}
\usepackage{xspace}

\usepackage{mathtools}          
\usepackage{amssymb}            
\usepackage{bm}                 
\usepackage{stmaryrd}           
\usepackage{accents}            
\usepackage{xfrac}              
\usepackage{booktabs}           
\usepackage{textcomp}           

\usepackage[foot]{amsaddr}
\usepackage[svgnames]{xcolor}
\usepackage{todonotes}

\calclayout

\newcommand{\doi}[1]{\textsc{doi}: \href{http://dx.doi.org/#1}{\nolinkurl{#1}}}

\DeclareMathAlphabet{\mathbsl}{OT1}{cmr}{bx}{sl}  

\newcommand{\mathsc}[1]{\text{{\rmfamily\scshape #1}}}

\newcommand {\V}  [1] {\bm{#1}}                        
\newcommand {\T}  [1] {\bm{#1}}                        
\newcommand {\TC} [1] {\mathbf{#1}}                    
\newcommand {\NM} [1] {\underaccent{\bar}{#1}}         

\newcommand   {\xa} {\hspace{-0.285em}}

\newcommand   {\xo} {\hspace{ 0.055em}}

\newcommand   {\jmpi} [1] {\llparenthesis\xo #1\xo\rrparenthesis}
\newcommand   {\jmpo} [1] {\left\llbracket #1\right\rrbracket}

\newcommand   {\avg}  [1] {\left\{\xa\left\{ #1 \right\}\xa\right\}}

\renewcommand {\d}            {\partial}
\newcommand   {\D}            {\mathrm{d}\mspace{1.0mu}}

\newcommand   {\transpose}[1] {#1^{\textsc{t}}}

\newcommand   {\Ne}           {N_{\mathsc{e}}}

\newcommand   {\CFL}          {\mathit{CFL}}
\newcommand   {\CFLc}         {\CFL^{\ast}}
\newcommand   {\CFLcLmb}      {\CFLc_{\lambda}}

\begin{document}

\title%
[
SDC for incompressible flow with variable viscosity
]
{
A spectral deferred correction method for incompressible flow with variable viscosity
}
\author{J\"org Stiller}
\address{%
  TU Dresden, Institute of Fluid Mechanics, 01062 Dresden, Germany}
\email{joerg.stiller@tu-dresden.de}

\begin{abstract}
This paper presents a semi-implicit spectral deferred correction (SDC) method 
  for incompressible Navier-Stokes problems with variable viscosity and 
  time-dependent boundary conditions.
The proposed method integrates elements of velocity- and pressure-correction 
  schemes, which yields a simpler pressure handling and a smaller splitting
  error than
  the SDPC method of Minion\,\&\,Saye (J.\ Comput.\ Phys.\ 375:\,797--822, 2018).
Combined with the discontinuous Galerkin spectral-element method for spatial
  discretization it can in theory reach arbitrary order  of accuracy in time 
  and space.
Numerical experiments in three space dimensions demonstrate up to order 12 
  in time and 17 in space for constant as well as 
  varying, solution-dependent viscosity.
Compared to SDPC the present method yields a substantial improvement of
  accuracy and robustness against order reduction caused by 
  time-dependent boundary conditions.

\end{abstract}

\keywords{%
Incompressible flow,
Variable viscosity,
High-order time integration,
Discontinuous Galerkin method.}

\maketitle


\section{Introduction}
\label{sec:intro}

High-order discretization methods are gaining interest in
  fluid mechanics 
  \cite{Schaal2015a_SE,Tavelli2016a_TI,Buscariolo2019a_SE,Fehn2019a_SE}, 
  solid mechanics 
  \cite{Delcourte2015a_TI,Tavelli2018a_TI}, 
  electrodynamics 
  \cite{Chen2013a_SE,Chan2017a_SE} and
  other areas of computational science 
  that are governed by partial differential equations, 
  such as meteorology and climate research 
  \cite{Marras2015a_SE}.
This development is driven by the expectation 
  of achieving a superior algorithmic efficiency
  which enables high-fidelity simulations at a scale 
  beyond the reach of low-order methods. 
In simulations of processes evolving in space-time, 
  the accuracy of spatial and temporal approximations must
  be tuned to each other.
The natural and only scalable way for achieving this is to match
  the convergence rates in space and time.
However, 
  examining recent work on high-order methods in 
  computational fluid dynamics (CFD) reveals that 
  the spatial order ranges typically from 4 to 16, 
  whereas the temporal order rarely exceeds 3.
This discrepancy constitutes the principal motivation of the present work,
  which adopts the spectral deferred correction (SDC) method to reach
  arbitrary temporal convergence rates for incompressible flows with variable
  viscosity.

Before reviewing the state of the art in high-order time integration,
  the implications of high-order space discretization shall be briefly 
  recapitulated.
In CFD, element based Galerkin methods with piecewise polynomial expansions
  represent the prevalent approach apart from spectral, finite-difference and
  isogeometric methods.
For an introduction and a comprehensive overview the reader is 
  referred to, e.g.,
  \cite{Deville2002a_SE,Karniadakis2005a_SE,Canuto2007a_SE,Hesthaven2008a_SE}.
The numerical properties of these methods have been thoroughly studied at
  the hand of convection, diffusion and wave problems as well as combinations
  thereof.
\citet{Ainsworth2009a_SE} and \citet{Gassner2011a_SE} 
  showed that high-order spectral element approximations achieve
  far lower dispersion and dissipation errors than
  second-order finite-element or finite-volume methods using a 
  comparable mesh spacing.
This property yields a tremendous advantage in marginally resolved
  simulations of turbulent flows
  \cite{Beck2014a_SE}.
On the other hand, the condition of the discrete operators worsens 
  when increasing the degree of the expansion basis.
Denoting the element size with $h$
  and the polynomial degree with $P$,
  the largest eigenvalues grow asymptotically as
  ${\lambda_{\mathsc c}\sim P^2 / h}$ 
  for the convection operator and
  ${\lambda_{\mathsc d} \sim P^4 / h^2}$ 
  for the diffusion operator,
  see \cite[Ch.~7.3]{Canuto2011a_SE}.
For explicit time integration schemes 
  these estimates imply stability restrictions of the form
  ${\Delta t_{\mathsc c} \sim h / P^2}$ and
  ${\Delta t_{\mathsc d} \sim h^2 / P^4}$,  
  respectively.
This corresponds to a reduction of the admissible time step
  by a factor of $1/P$ for convection and $1/P^2$ for diffusion
  in comparison to low-order finite-element or finite-volume methods.
These stability issues lead to a preference of implicit 
  time integration methods, especially for diffusion.
With convection, however, nonlinearity complicates
  implicit methods and renders semi-implicit or even fully explicit 
  approaches attractive.

Time integration methods applied in computational fluid dynamics
  cover a wide range of approaches, including
  multistep, 
  Runge-Kutta (RK), 
  Rosenbrock
  and 
  extrapolation methods
  \cite{Hairer1993a_TI,%
        Hairer1996a_TI,%
        Gottlieb2001a_TI} 
  as well as
  variational methods
  \cite{Thomee2006a_TI}.
%
%
Early work on high-order space discretization for incompressible flows
  advocated semi-implicit multistep methods 
  \cite{Orszag1986a_TI,Karniadakis1991a_TI} 
  based on the projection method introduced by 
  \citet{Chorin1968a_TI}.
These methods were investigated and generalized in numerous follow-up
  studies, and gained considerable popularity, mainly because of their 
  simplicity and low cost per time step, e.g.,
  \cite{Guermond2006a_TI,%
        Leriche2006a_TI,%
        Ferrer2014a_TI,%
        Fehn2017a_TI}.
Their advantage is offset, however, by
  harsh stability restrictions for the explicit part and
  by the loss of A-stability of the implicit part 
  for convergence orders greater than two
  (\emph{second Dahlquist barrier} \cite{Hairer1996a_TI}).
Accordingly, the vast majority of studies based on multistep methods
  uses order two in time, whereas the spatial order ranges from 4 to
  well above 10. 
As a notable exception \citet{Klein2015a_TI} used a fully implicit 
  SIMPLE method based on backward differentiation formulas (BDF) up 
  to order 4.

%
%
In contrast to linear multistep methods, implicit Runge-Kutta methods 
  can be constructed to reach high convergence orders along with excellent
  stability properties.
For flow problems, however, convection or solution-dependent viscosity 
  result in nonlinear equations which need to be linearized and solved on
  every stage.
In order to reduce complexity, all studies known to the author
  used diagonally implicit Runge-Kutta (DIRK) methods.
\citet{Uranga2010a_TI}
  applied a third order DIRK for large-eddy simulations of
  transitional flow wings in conjunction with 
  Newton's method and preconditioned conjugate gradients 
  for solving the nonlinear equation systems on each stage.
Rosenbrock-type methods are build on the Jacobian of the 
  right-hand side (RHS) and, thus, 
  achieve linearization more directly.
\citet{John2006a_TI} compared Rosenbrock methods of order 3
  for incompressible Navier-Stokes problems with fractional step
  methods and showed their competitiveness, especially in terms 
  of accuracy and robustness.
More recently, 
  \citet{Bassi2015a_TI} and
  \citet{Noventa2016a_TI}
  applied Rosenbrock methods of orders up to 6 
  combined with discontinuous Galerkin methods in space
  to compressible and incompressible flows past airfoils
  and other configurations.
As an alternative to the fully implicit approach, 
  implicit-explicit (IMEX) methods combine 
  implicit RK for the stiff (and often linear) part of the RHS 
  with
  explicit RK schemes for the nonlinear part.
Following
  \citet{Ascher1997a_TI}, who presented IMEX RK methods up to order 3,
  \citet{Kennedy2003a_TI} developed methods up to order 5 with embedded
  lower order schemes for error estimation.
More recently,
  \citet{Cavaglieri2015a_TI} devised IMEX RK methods with reduced memory
  requirements, while
  \citet{Boscarino2016a_TI} consider the extension to problems which do not 
  allow for a sharp separation of the stiff RHS part.%
So far, applications of IMEX RK to flow simulations featuring 
  high-order spatial approximations seem to be rare and confined
  to the compressible case, e.g., 
  \cite{Persson2011a_TI,Higueras2014a_TI,Gardner2018a_TI}.
Moreover, the increasing number of matching conditions \cite{Kennedy2003a_TI}
  complicates the construction of higher order IMEX RK and methods of order 
  higher than five are not available to the knowledge of the author.
A further issue arising in the construction of high-order RK methods 
  is the order reduction phenomenon, which can be triggered, e.g., 
  by
  stiff source terms \cite{Carpenter1995a_TI,Hairer1996a_TI} or 
  time-dependent boundary conditions \cite{Kennedy2003a_TI}.
To alleviate this issue several approaches have been proposed 
  \cite{Alonso-Mallo2002a_TI,Alonso-Mallo2005a_TI,Rosales2017a_TI}, 
  but no solution is known for complex problems such as the 
  Navier-Stokes equations.
One possibility would be to use methods possessing a high stage order
  \cite{Burrage1990a_TI,Hairer1996a_TI}.
But unfortunately,
  DIRK and Rosenbrock methods are limited to stage order 
  two by construction.

%
%
Extrapolation and deferred correction methods employ low-order 
  time-integration schemes within an iterative framework to achieve
  convergence of (in principle) arbitrary high order, 
  see
  \cite{Constantinescu2010a_TI} and
  \cite{Dutt2000a_TI}, respectively.
The spectral deferred correction (SDC) method was developed by 
  \citet{Dutt2000a_TI} for solving the Cauchy problem for ordinary 
  differential equations and extended by
  \citet{Kress2002a_TI} to initial boundary value problems.
\citet{Minion2003b_TI} generalized the originally either implicit or 
  explicit approach by proposing a semi-implicit SDC method.
The basic idea of SDC is to convert the differential evolution 
  problem into a Picard integral equation which is solved by a deferred 
  correction procedure, driven by a lower order marching scheme.
In this procedure, the lower order scheme sweeps repeatedly 
  through subintervals defined by a set of collocation points,
  which also serve for Lagrange interpolation and integration.
Choosing these points from a Gauss-type quadrature yields an SDC
  method that converges toward the solution of the corresponding
  implicit Gauss collocation method 
  \cite{Hairer1996a_TI}.
Applying a first order corrector such as implicit Euler,
  each sweep ideally increases the order by one, until 
  reaching the maximum depending on the chosen set of
  collocation points.
For elevating the order by more than one per sweep, correctors based 
  on multistep and RK methods have been considered, e.g.,  
  \cite{Layton2008a_TI,%
        Christlieb2009a_TI,%
        Christlieb2009b_TI,%
        Christlieb2011a_TI}.
\citet{Christlieb2009b_TI} showed, however, that this approach 
  imposes smoothness conditions on the error vector that are 
  difficult, if possible, to meet with nonuniform (Gauss-type) 
  points.
As a remedy the authors proposed methods based on uniform points
  and embedded high-order RK integrators.
These methods, indeed, achieve the expected higher order improvement 
  per sweep, but are limited to roughly half the order attained by
  a Gauss-type method using the same number of points.
Other approaches to accelerate the SDC method include
	the application of high-order schemes for computing the initial 
	approximation \cite{Layton2007a_TI},
	preconditioning using Krylov subspace methods \cite{Huang2006a_TI}
	or optimized DIRK-type sweeps \cite{Weiser2015a_TI} and
	multi-level SDC \cite{Speck2015a_TI,Speck2016a_TI}.
Like other high-order methods SDC is susceptible to order reduction,
  which can be caused by stiff source terms 
  \cite{Huang2006a_TI,Christlieb2015a_TI}
  or boundary conditions
  \cite{Minion2018a_TI}.
However, in contrast to RK methods, the phenomenon manifests itself
  in a slower convergence of the correction sweeps
  rather than in a reduction of the order of the final solution.  
In spite of its capability to reach arbitrary high orders and 
  straightforward extension to parallel-in-time methods boosting 
  the efficiency on high-performance computers
  \cite{Minion2015a_TI,Bolten2017a_TI},
  SDC has been rarely applied to fluid dynamics problems.
Moreover, most studies were confined to simple configurations
  with periodic boundaries
  \cite{Minion2003a_TI,%
        Minion2004a_TI,%
        Almgren2013a_TI,%
        Speck2015a_TI}.
Only recently,
  \citet{Minion2018a_TI}
  proposed a semi-implicit SDC method based on a first-order projection 
  scheme for 2D incompressible flows.
The implicit part of their method is accelerated by DIRK sweeps
  as proposed by \citet{Weiser2015a_TI}.
Nevertheless, it suffered from severe order reduction when applied 
   with time-dependent Dirichlet conditions.

%
%
Variational methods resemble SDC in harnessing piecewise polynomial expansions
  in the time direction.
Unlike the latter, they achieve discretization by application of a variational
  principle, predominantly the discontinuous Galerkin (DG) method.
Recent applications of DG in time include
  incompressible flows \cite{Tavelli2016a_TI,Ahmed2017a_TI},
  elasticity \cite{Tavelli2018a_TI} 
  and hyperbolic conservation laws \cite{Friedrich2019a_TI}.
Although inherently implicit, the method allows to incorporate semi-implicit
  strategies similar to SDC as proposed e.g.\ in \cite{Tavelli2016a_TI}.
Like SDC, DG methods can be based on Lagrange polynomials constructed from
  Gauss points.
However, assuming that ${Q\!+\!1}$ points are used, 
  DG methods generally converge with order ${Q\!+\!1}$, 
  whereas SDC methods reach order
  ${2Q}$ with Gauss-Lobatto-Legendre (GLL) and
  ${2Q\!+\!2}$ with Gauss-Legendre points
  \cite{Causley2019a_TI}.

%
%
The goal of this study is to develop an SDC method for incompressible 
  Navier-Stokes problems which, in combination with the 
  DG spectral element method for spatial discretization, 
  is capable to reach arbitrary high order in time and space.
As the backbone of the new method, a semi-implicit correction scheme
  is devised which yields a simpler and more robust pressure handling than 
  the SDPC method proposed by \citet{Minion2018a_TI}.
In contrast to the latter, which is based on a pure pressure correction scheme,
  the proposed method combines ideas of velocity and pressure correction to 
  reduce the splitting error and to extend the approach to variable viscosity.
Additionally, it penalizes the divergence inside and jumps across the elements
  to improve the continuity of the approximate velocity. 
The effect of these measures is confirmed in numerical studies which reveal 
  a substantial improvement over SDPC, especially in the case of time-dependent
  boundary conditions.
Moreover, the proposed SDC method is shown to work equally well with
  variable and even solution-dependent viscosity.

%
%
The remainder of the paper is organized as follows: 
Section \ref{sec:navier-stokes} summarizes the incompressible Navier-Stokes
  equations with variable viscosity.
Section \ref{sec:sdc} reviews the spectral deferred correction method and 
  extends it to the flow problem.
Section \ref{sec:spatial discretization} presents the spatial discretization
  followed by a compilation and discussion of numerical results in section
  \ref{sec:numerical experiments}.
Section \ref{sec:conclusions} concludes the paper.


\section{The Navier-Stokes equations with variable viscosity}
\label{sec:navier-stokes}

This paper considers incompressible flows with constant density and variable
  viscosity in a simply connected domain ${\Omega \in \mathbb R^3}$.
The velocity ${\V v(\V x, t)}$ satisfies the momentum (Navier-Stokes) and
  continuity equations
  \begin{gather}
    \label{eq:momentum}
    \d_t \V v + \nabla \cdot \V v \V v + \nabla p
    = \nabla \cdot \T\tau + \V f
    \,,
  \\
    \label{eq:continuity}
    \nabla \cdot \V v = 0
  \end{gather}
  in $\Omega$,
  where 
  $p$ represents the pressure and
  \begin{equation}
    \label{eq:stress tensor}
    \T\tau = \nu [\nabla \V v + \transpose{(\nabla \V v)}] 
  \end{equation}
  the viscous stress tensor, 
  both divided by density;
  ${\nu(\V x, t, \V v)}$ is the kinematic viscosity and 
  ${\V f(\V x, t)}$ an explicitly defined forcing term.
The flow problem is closed be stating 
  initial and boundary conditions
  \begin{alignat}{2}
  \label{eq:ic}
  \V v(\V x, 0) &= \V v_0(\V x)
  &\quad& \V x \in  \Omega
  \,,
  \\
  \label{eq:bc}
  \V v(\V x, t) &= \V v_{\mathrm b}(\V x, t)
  &\quad& \V x \in \d\Omega
  \,.
  \end{alignat}
For continuity, 
  $\V v_0$ must be divergence free and
  $\V v_{\mathrm b}$ satisfy the compatibility condition
  \begin{equation}
    \int_{\d\Omega} \V n \cdot \V v_{\mathrm b} \D\Gamma = 0
    \,.
  \end{equation}

Assuming a constant viscosity and using the identity
  ${\nabla^2 \V v  + \nabla \times \nabla \times \V v 
   = \nabla\nabla \cdot \V v }$
  the viscous term in the momentum equation \eqref{eq:momentum} 
  can be rewritten in the following forms
  \begin{subequations}
    \label{eq:viscous term}
    \begin{alignat}{2}
      \label{eq:viscous term:native}
      \nabla \cdot \T\tau
      & = \phantom{-}
           \nu [\nabla^2 \V v + \nabla\nabla \cdot \V v] 
      &\qquad
      & \text{(native)}
    \\
      \label{eq:viscous term:laplacian}
      & = \phantom{-}
           \nu \nabla^2 \V v                           
      &
      & \text{(laplacian)}
    \\
      \label{eq:viscous term:rotational}
      & = -\nu \nabla \times \nabla \times \V v^{\phantom{2}}
      &
      & \text{(rotational)}\,. 
    \end{alignat}
  \end{subequations}
For variable $\nu$ they can be generalized to
  \begin{equation}
     \label{eq:viscous term:general}
     \nabla \cdot \T\tau
     = \nabla \cdot \nu [ \nabla \V v  
                        + \transpose{(\nabla \V v)} 
                        - \chi\TC I \nabla \cdot \V v
                        ]
     \,,
  \end{equation}
  where 
  ${\chi = 0}$ corresponds to the native, 
  ${\chi = 1}$ the laplacian, and
  ${\chi = 2}$ the rotational form, respectively.
These forms are equivalent when applied to solenoidal vector fields, 
  but not for approximate solutions that are not divergence free.
This needs to be considered in the discrete case.


\section{Spectral deferred correction method}
\label{sec:sdc}


\subsection{General approach}


This section briefly reviews the SDC method based on the model problem
  \begin{equation}
    \label{eq:model problem}
    \D_t\V v = \V F(t, \V v(t))\,
  \end{equation}
  with
  ${t \in (t_0,T]}$,
  ${T = t_0 + \Delta t}$
  and initial condition ${\V v(t_0) = \V v_0}$.


\subsubsection{Preliminaries}

As a prerequisite for developing the method,
  the time interval is divided into subintervals 
  ${\{(t_{i-1}, t_i)\}_{i=1}^M}$
  such that 
  ${t_0 = t^n}$
  and 
  ${t_M = t^{n+1}}$.
Further, let
  ${\Delta t_i = t_i - t_{i-1}}$
  denote the length of the $i$-th subinterval,
  ${\V v_i \simeq \V v(t_i)}$ 
  the discrete solution at time $t_i$
  and
  $\V v_i^k$ the \mbox{$k$-th} approximation of $\V v_i$.
In addition to this, the intermediate times $t_i$ serve
  as collocation points for Lagrange interpolation and
  as quadrature points for numerical integration.
Depending on the underlying quadrature rule, one or both endpoints
  may be dropped, see e.g.\ \cite{Layton2005a_TI,Causley2019a_TI}.
The following description is based on the Gauss-Lobatto-Legendre (GLL)
  rule and, hence, includes both endpoints for interpolation and
  integration. 
Accordingly, 
  $\{t_i\}_{i=0}^M$ represent the GLL points scaled to ${[t_0,T]}$,
  ${\NM{\V v} = [ \V v_i ]}$ 
  the discrete solution vector,
  ${\NM{\V v}^k = [\V v_i^k]}$,
  and
  ${\mathcal I \NM{\V v}^k(t)}$ 
  the corresponding Lagrange interpolant at time $t$. 


\subsubsection{Predictor}

The initial approximation $\NM{\V v}^0$ is obtained by performing a 
  predictor sweep of the form
  \begin{equation}
    \label{eq:predictor}
    \V v^0_i = \V v^0_{i-1} + \V H_i(\NM{\V v}^0) 
    \, , \quad i=1, \dots, M 
    \,,
  \end{equation}
  where
  ${\V v^0_0 = \V v_0}$
  and
	${\V H_i(\NM{\V v})}$ is an approximation of 
	${\int_{t_{i-1}}^{t_i}\!\V F \D t}$.
Using, for example, a combination of forward and backward Euler rules
  based on the decomposition
  ${\V F = \V F^{\mathsc{im}} + \V F^{\mathsc{ex}}}$
  yields
  an IMEX Euler predictor
  with
  \begin{equation}
    \V H_i(\NM{\V v}) 
    = \Delta t_i \big[ \V F^{\mathsc{im}}(t_{i}  , \V v_{i})
                     + \V F^{\mathsc{ex}}(t_{i-1}, \V v_{i-1}) 
                 \big]
    \,.
  \end{equation}
Alternatively, the predictor can be constructed from higher order 
  time integration schemes such as RK or multistep methods
  \cite{Layton2007a_TI}.
This approach may give an advantage by providing more accurate starting values,
  but will not be investigated in frame of the present study.


\subsubsection{Corrector}

The goal of the corrector is to remove the error from a given approximation
  $\NM{\V v}^k$.
For deriving the correction equation,
  the error function is defined as
  \begin{equation}
    \label{eq:error function}
    \V\delta^k(t) = \V v(t) - \mathcal I \NM{\V v}^k(t) 
    \,.
  \end{equation}
Further, the residual function is introduced by
  \begin{equation}
    \label{eq:residual function}
    \V \varepsilon^k(t)
    = \V v_0 
    + \int_{t_0}^t \V F(\tau,\mathcal I \NM{\V v}^k(\tau))\,\D\tau 
    - \mathcal I \NM{\V v}^k(t)
    \,.
  \end{equation}
Differentiating and subtracting
  \eqref{eq:error function} and
  \eqref{eq:residual function}
  yields
  \begin{equation}
    \D_t (\V\delta^k - \V\varepsilon^k)
      = \D_t \V v(t) - \V F(t, \mathcal I \NM{\V v}^k(t))
      \,.
  \end{equation}
This equation can be rearranged using 
  \eqref{eq:model problem} and 
  \eqref{eq:error function}
  to give the error equation
  \begin{equation}
    \label{eq:error equation}
    \D_t (\V\delta^k - \V\varepsilon^k)
      = \V F \big(t, \mathcal I \NM{\V v}^k(t)\!+\!\V\delta^k(t) \big) 
      - \V F \big(t, \mathcal I \NM{\V v}^k(t)\big)
      \,,
  \end{equation}
  which is supplemented with the initial condition 
  ${\V \delta^k(t_0) = 0}$.

The error equation is solved numerically by means of a time integration
  scheme sweeping through the subintervals. 
For example, application of IMEX Euler yields
  \begin{equation}
    \label{eq:approx error}
    \V\delta^k_i = \V\delta^k_{i-1} 
                 + \V \varepsilon^k(t_i)
                 - \V \varepsilon^k(t_{i-1})
                 + \V H_i(\NM{\V v}^k + \NM{\V \delta}^k) 
                 - \V H_i(\NM{\V v}^k) 
  \end{equation}
  for
  ${i=1, \dots, M}$,
  where 
  $\V\delta^k_i$ 
  represents the approximation of
  $\V\delta^k(t_i)$.
Substituting
  \eqref{eq:residual function}
  for 
  ${\V \varepsilon^k(t_{\ast})}$, 
  eliminating
  $\V \delta^k_{\ast}$
  by means of \eqref{eq:error function}
  and
  defining 
  the new approximate solution by
  ${\V v^{k+1}_i \!= \V v^k_i + \V\delta^k_i}$
  finally gives the 
  update equation
  \begin{equation}
    \label{eq:corrector}
    \V v^{k+1}_i
      = \V v^{k+1}_{i-1}
      + \V H_i(\NM{\V v}^{k+1}) 
      - \V H_i(\NM{\V v}^k) 
      + \int_{t_{i-1}}^{t_i} \!
          \V F(\tau,\mathcal I \NM{\V v}^k(\tau))\,\D\tau
    \,.
  \end{equation}   
The last term in \eqref{eq:corrector} is usually approximated by
  replacing the integrand by its Lagrange interpolant, i.e.
  \begin{equation}
    \int_{t_{i-1}}^{t_i} \!
      \V F(\tau,\mathcal I \NM{\V v}^k(\tau))\,\D\tau
    \approx
      \int_{t_{i-1}}^{t_i} \! \mathcal I \NM{\V F}^k(\tau)\,\D\tau
    \eqqcolon
    \V S^k_i
    \,, 
  \end{equation}
  with 
  ${\NM{\V F}^k = [\V F(t_i,\V v^k_i)]}$.
The approximate integral can be expressed in terms of a quadrature formula,
  \begin{equation}
    \label{eq:S}
    \V S^k_i 
      = \Delta t \sum_{j=0}^M w^{\mathsc s}_{i,j} \V F(t_j,\V v^k_j)
    \,,
  \end{equation}
  where
  $w^{\mathsc s}_{i,j}$
  equals the integral of the 
  $j$-th interpolation polynomial
  over subinterval ${(t_{i-1},t_i)}$,
  normalized with $\Delta t$.
As a consequence, the sum
  $\smash{\sum_i w^{\mathsc s}_{i,j}}$
  recovers the weights of the underlying quadrature rule
  and, hence,
  $\NM{\V v}^k$ converges to the solution
  of the related collocation method.

The sketched SDC method attains 
  order ${2M+1}$ at the final time $T$ of a single interval
  and $2M$ when repeated for stepping through a sequence of 
  multiple intervals.
Using a first order corrector as sketched above, 
  every sweep, ideally, elevates the order
  by one, until reaching the maximum order
  \cite{Dutt2000a_TI}.
However, stiff terms and boundary conditions may affect convergence
  such that more iterations are required to attain the optimal order.


\subsection{Application to incompressible flow}


\subsubsection{Considerations}


The generalization of the SDC approach 
  to the incompressible Navier-Stokes problem (\ref{eq:momentum}--\ref{eq:bc})
  follows a similar approach as outlined above.
Starting from an identical partitioning of a given time interval,
  the semi-discrete solution is denoted by
  ${\V v_i(\V x) \simeq \V v(\V x, t_i)}$ 
  for the velocity 
  and
  ${p_i(\V x) \simeq p(\V x, t_i)}$ 
  for the pressure.
Similarly, ${\V v_i^k(\V x)}$ and $p_i^k(\V x)$
represent the corresponding approximations after
$k$ correction sweeps.
Before proceeding it is important to note 
  the following differences between 
  the flow problem and 
  the model problem
  \eqref{eq:model problem}:
Although the momentum balance \eqref{eq:momentum} resembles an evolution
  equation for the velocity, it involves an additional variable in terms
  of the pressure. 
Moreover, the velocity is required to satisfy the continuity equation
  \eqref{eq:continuity}, which lacks a time derivative and, hence, looks
  like an algebraic constraint from perspective of time integration.
Finally, the flow equations are subject to boundary conditions 
  \eqref{eq:bc} that may depend on time themselves.


\subsubsection{Predictor}


The complex nature of the flow problem complicates the
  construction of the predictor \eqref{eq:predictor}.
Instead of defining the operator $\V H_i$ directly it is more 
  appropriate to derive its structure from a single time step 
  across some subinterval ${(t_{i-1},t_i)}$.
In analogy to the model problem, on could use the IMEX Euler method 
  for incompressible flow, i.e. 
  \begin{gather}
    \label{eq:euler:momentum}
    \frac{\V v_{i} - \V v_{i-1}}{\Delta t_i}
    + \nabla \cdot (\V v \V v)_{i-1} + \nabla p_i 
    = \nabla \cdot \T \tau_i 
    + \V f_i
    \,,
  \\
    \label{eq:euler:continuity}
    \nabla \cdot \V v_i = 0
    \,.
  \end{gather}
While this scheme looks reasonably simple and elegant,
  it yields a coupled system for $\V v_i$ and $p_i$,
  which renders the solution costly, 
  especially in view of the
  pertinent stability restrictions  
  and low accuracy.
Therefore, it seems attractive turning to projection schemes
  that decouple continuity from the momentum balance.
These schemes employ some approximation of pressure $p_i$ in the
  (incomplete) momentum step and 
  achieve continuity by performing a separate projection step 
  \cite{Glowinski2003a_TI,Guermond2006a_TI,Minion2018a_TI}.
Moreover, the latter yields a Poisson equation for correcting the pressure.

Depending on the order of the substeps two classes of projection schemes
  can be distinguished: pressure-correction and velocity-correction methods.
Pressure-correction methods were introduced by \citet{Chorin1968a_TI}.
They first solve an implicit diffusion problem for each velocity component, 
  including approximations for convection and pressure terms, and then 
  project the provisional velocity to a divergence-free field.
In comparison to IMEX Euler, the splitting leads to an additional
  error caused by the violation of tangential boundary conditions 
  in the projection step.
However, several approaches exist for controlling the splitting error and 
  retaining first order convergence \cite{Guermond2006a_TI}.
\citet{Minion2018a_TI} investigated different variants of the pressure-correction
  scheme as a basis for their SDPC method.

Velocity-correction methods were introduced by \citet{Orszag1986a_TI}
  and further extended, e.g. in \cite{Karniadakis1991a_TI,Guermond2003b_TI}.
As a common feature, these methods perform the projection step before
  solving the diffusion problem.
With constant viscosity this approach preserves continuity.
However, similar to the pressure-correction method, it introduces 
  a splitting error due to inaccurate pressure boundary conditions 
  in the projection step.
Using the rotational form of the velocity-correction method
  mitigates this error and recovers optimal convergence
  \cite{Guermond2003b_TI}.

Several authors adapted the pressure-correction method to simulate 
  flows with variable viscosity
  \cite{Fai2014a_TI,Niemann2018a,Deteix2018a_TI},
  whereas the author is not aware of corresponding extensions of the
  velocity-correction method.
On the other hand, the pressure-correction method implies a rather 
  complicated handling of the pressure when applied as a base method
  for SDC \cite{Minion2018a_TI}.
The method proposed in the following combines the advantages of both 
  approaches:
It starts with a velocity-correction step and concludes with a projection
  like the pressure-correction method.


For stating the base time-integration method, the momentum equation is
  rewritten in the form
  \begin{equation}
    \d_t \V v = \V F(\V x, t, \V v, p)
    \,,
  \end{equation}
  where
  \begin{equation}
    \label{eq:F}
    \V F = \V F_c(\V v) 
         + \V F_d(\nu, \V v) 
         + \V F_p(p) 
         + \V f(\V x, t)
  \end{equation}
  with
  \begin{alignat}{3}
    \label{eq:F_c}
    & \V F_c    && = &\,-& \nabla \cdot \V v \V v                        \\
    \label{eq:F_d}
    & \V F_d    && = &   & \V F_{d1} + \V F_{d2} + \V F_{d3}             \\
    \label{eq:F_d1}
    & \V F_{d1} && = &   & \nabla \cdot \nu \nabla \V v                  \\
    \label{eq:F_d2}
    & \V F_{d2} && = &   & \nabla \cdot \nu \transpose{(\nabla \V v)}    \\
    \label{eq:F_d3}
    & \V F_{d3} && = &\,-& \chi \nabla (\nu \nabla \cdot \V v)           \\
    \label{eq:F_p}
    & \V F_{p}  && = &\,-& \nabla p
    \,.
  \end{alignat}
These definitions trivially extend to the semi-discrete solution.
For brevity the arguments are omitted whenever possible, as for example in
$\V F_{d,i}^k = \V F_d(\nu_i^k, \V v_i^k)$.

The predictor is then defined as follows:

  \begin{alignat}{2}
    \label{eq:predictor:1}
    &
    \frac{\V v'^{\,0}_i - \V v^0_{i-1}}{\Delta t_i}
      && = (\V F_{c} + \V F_{d})^0_{i-1}
         + \V f_{i}
    \,,
  \\
    \label{eq:predictor:2}
    &
    \frac{\V v''^{\,0}_i - \V v'^{\,0}_i}{\Delta t_i} 
      && = \V F_{p}(p''^{\,0}_i), \quad
    \nabla \cdot \V v''^{\,0}_i = 0, \quad
    \V n \cdot \V v''^{\,0}_i|_{\d\Omega} = \V n \cdot \V v_{\textrm b}(t_i)
    \,,
  \\
    \label{eq:predictor:3}
    &
    \frac{\V v'''^{\,0}_i - \V v''^{\,0}_i}{\Delta t_i} 
      && = \V F_{d1}(\nu^0_{i-1},\V v'''^{\,0}_i)
         - (\V F_{d1} + c_{\chi}\V F_{d3})^0_{i-1},
    \quad
    \V v'''^{\,0}_i|_{\d\Omega} = \V v_{\textrm b}(t_i)
    \,,
  \\
    \label{eq:predictor:4}
    &
    \frac{\V v^0_i - \V v'''^{\,0}_i}{\Delta t_i} 
      && =   \V F_p(p^0_i - p''^{\,0}_i),
    \quad
    \nabla \cdot \V v^0_i = 0,
    \quad
    \V n \cdot \V v^0_i|_{\d\Omega} = \V n \cdot \V v_{\textrm b}(t_i)
    \,,
  \end{alignat}
  for ${i = 1, \dots, M}$.
The first three substeps comprise a velocity-correction method.
In particular, \eqref{eq:predictor:1} represents an incomplete Euler step
  using the forward rule for convection and diffusion,
  backward rule for the forcing term 
  and
  skipping the pressure part.
It is followed by the first projection \eqref{eq:predictor:2}
  and the viscous correction \eqref{eq:predictor:3}.
Note that the latter drops 
the diffusion term $\V F_{d1,i-1}^0$ and
a scaled part of the divergence contribution $\V F_{d3,i-1}^0$ 
introduced in the first substep.
The scaling factor is set to ${c_{\chi} = \sfrac{1}{2}}$ for ${\chi=2}$
and zero otherwise.
For constant viscosity
  (\ref{eq:predictor:1} -- \ref{eq:predictor:3})
  reproduce
  with ${\chi=1}$ the standard and
  with ${\chi=2}$ the rotational velocity-correction method
  as defined in \cite{Guermond2003b_TI}.
In the case of variable viscosity,
  the diffusion step produces a divergence error 
  of the order ${O(\Delta t_i)}$.
This error is removed by the final projection step \eqref{eq:predictor:4}. 
Alternatively, it can be tolerated as a part of the overall discretization 
  error, which will be considered as an option in the numerical experiments.

In contrast to pressure-correction, the proposed method requires no
  initial approximation of the pressure.
The intermediate pressure ${p''^{\,0}_i}$ and the final pressure ${p^0_i}$
  are obtained each by solving a Poisson problem 
  which follows from the corresponding projection step.
For example, taking the divergence and, respectively, 
  the normal projection of the first equation in \eqref{eq:predictor:2} 
  leads to
  \begin{gather}
    \label{eq:pressure problem}
    \nabla^2 p''^{\,0}_i = \frac{1}{\Delta t_i} \nabla \cdot \V v'^{\,0}_i
    \,, \quad
    \V n \cdot \nabla p''^{\,0}_i|_{\d\Omega} = 
      \frac{1}{\Delta t_i}\,
      \V n \cdot (\V v'^{\,0}_i|_{\d\Omega} - \V v_{\textrm b}(t_i))
    \,.
  \end{gather}
Similarly, \eqref{eq:predictor:4} yields a Poisson equation 
  for ${p^0_i - p''^{\,0}_i}$
  with homogeneous Neumann conditions.


\subsubsection{Corrector}
\label{sec:sdc:flow:corrector}

Apart from the additional low- and high-order contributions, 
  the corrector resembles the predictor.
The substeps are
  \begin{alignat}{2}
    \label{eq:corrector:ansatz:1}
    &
    \frac{\V v'^{\,k+1}_i - \V v^{k+1}_{i-1}}{\Delta t_i}
      && = ( \V F_{c} 
           + \V F_{d2}
           + \V F_{d3}
           )_{i-1}^{k+1}
           + \V F_{d1}(\nu^{k+1}_{i-1},\V v^{k}_i)
           - \frac{\V H^{k}_i}{\Delta t_i}
           + \frac{\V S^{k}_i}{\Delta t_i}
    \,,
  \\
    \label{eq:corrector:ansatz:2}
    &
    \frac{\V v''^{\,k+1}_i - \V v'^{\,k+1}_i}{\Delta t_i} 
      && = \V F_{p}(\tilde p''^{\,k+1}_i), \quad
    \nabla \cdot \V v''^{\,k+1}_i = 0, \quad
    \V n \cdot \V v''^{\,k+1}_i|_{\d\Omega} = \V n \cdot \V v_{\textrm b}(t_i)
    \,,
  \\
    \label{eq:corrector:ansatz:3}
    &
    \frac{\V v'''^{\,k+1}_i - \V v''^{\,k+1}_i}{\Delta t_i} 
      && = \V F_{d1}(\nu^{k+1}_{i-1}, \V v'''^{\,k+1}_i - \V v^{\,k}_i)
         - c_{\chi}\V F_{d3,i-1}^{k+1},
    \quad
    \V v'''^{\,k}_i|_{\d\Omega} = \V v_{\textrm b}(t_i)
    \,,
  \\
    \label{eq:corrector:ansatz:4}
    &
    \frac{\V v^k_i - \V v'''^{\,k}_i}{\Delta t_i} 
      && = \V F_p ( \tilde p^k_i 
                   - \tilde p''^{\,k}_i),
    \quad
    \nabla \cdot \V v^k_i = 0,
    \quad
    \V n \cdot \V v^k_i|_{\d\Omega} = \V n \cdot \V v_{\textrm b}(t_i)
    \,,  
  \end{alignat}
  where 
  \begin{equation}
     \label{eq:S:ansatz}
     \V S^{k}_i 
        = \Delta t 
           \sum_{j=0}^M w^{\mathsc s}_{i,j} 
             \Big[ \V F_{c,i}^{k} 
                 + \V F_{d,i}^{k}
                 + \V f_{i}
             \Big]
           - \nabla P^{k}_i
  \end{equation}
  represents the subinterval integral similar to \eqref{eq:S}
  with the pressure part $P^{k}_i$ yet to be defined.
In contrast to the predictor, the corrector exploits the previous 
  approximation of $\V v_i$ to provide a more accurate starting value for
  the implicit diffusion term in the extrapolation step
  \eqref{eq:corrector:ansatz:1}.
As will become clear in a moment,
  the pressure computed in the projection steps 
  \eqref{eq:corrector:ansatz:2} and
  \eqref{eq:corrector:ansatz:4} 
  is in general not an approximation of $p_i$ 
  and, hence, marked by a tilde.

Adding and rearranging the equations for 
  $\V v'^{\,k+1}_i$, 
  $\V v''^{\,k+1}_i$, 
  $\V v'''^{\,k+1}_i$ 
  and 
  $\V v^{k+1}_i$
  gives
  \begin{equation}
    \label{eq:corrector:full-step:ansatz}
    \begin{aligned}
      \V v^{k+1}_i - \V v^{k+1}_{i-1}
        & = \Delta t_i
            \Big[ \V F_{c,i-1}
              + \V F_{d1}(\nu_{i-1}, \V v'''_i)
              + \V F_{d2,i-1}
              + (1 - c_{\chi}) \V F_{d3,i-1}
              + \V F_{p}(\tilde p_i) 
          \Big]^{k+1}
        \\
        & - \V H^{k}_i
          + \V S^{k}_i
    \end{aligned}
  \end{equation}
and leads to the following ansatz for the low-order contribution 
  \begin{gather}
     \label{eq:H:ansatz}
     \V H^{k}_i 
       = \Delta t_i 
            \Big[ \V F_{c,i-1}
                + \V F_{d1}(\nu_{i-1}, \V v'''_i)
                + \V F_{d2,i-1}
                + (1 - c_{\chi}) \V F_{d3,i-1}
                + \V F_{p}(p_i^{\mathsc h})
            \Big]^{k}
    \,.
  \end{gather}
Substituting the high- and low-order contributions
  (\ref{eq:S:ansatz},\,\ref{eq:H:ansatz})
  in
  \eqref{eq:corrector:full-step:ansatz} 
  and considering the converged case
  yields
  \begin{gather}
    \label{eq:converged case}
    \V v_i 
      - \V v_{i-1} 
      = \Delta t 
           \sum_{j=0}^M w^{\mathsc s}_{i,j} 
             \big[ \V F_{c,i} 
                 + \V F_{d,i}
                 + \V f_{i}
             \big]
      - \nabla\big[ P_i + \Delta t_i\,\tilde p_i \big]
      \,.
  \end{gather}
Comparing this result to the corresponding collocation formulation
  implies 
  \begin{equation}
    \label{eq:pressure term}
    P_i + \Delta t_i\,\tilde p_i
    = \Delta t \sum_{j=0}^M w^{\mathsc s}_{i,j} p_i
    \approx
      \int_{t_{i-1}}^{t_i} \!\! p\,\D t
    \,.
  \end{equation}
Since $p_i^{\mathsc h}$ vanishes from 
  (\ref{eq:converged case},\ref{eq:pressure term})
  and 
  $P_i$ is balanced by $\tilde p_i$,
  the choice of these quantities seems to have no effect 
  on the SDC method.
This conjecture is confirmed by preliminary studies exploring
  several approaches, including the evaluation of 
  ${P_i^{k} = \Delta t \sum_j w^{\mathsc s}_{i,j} p^{k}_i}$ 
  using the recomputed pressure obtained from
  \begin{align}
    \label{eq:p:consistent:pde}
    & \nabla^2 p^{k}_i 
        = \nabla \cdot ( \V F_{c} 
                       + \V F_{d} 
                       + \V f    )^{k}_{i}
    \,,
    \\
    \label{eq:p:consistent:bc}
    &  \d_{\V n} p^{k}_i
         = \V n \cdot [ ( \V F_{c} 
                      + \V F_{d}
                      + \V f     )^{k}_{i}
                      - \d_t \V v_{\textrm b, i} ]_{\d\Omega}
    \,.
  \end{align}
Consequently, all studies in this work were performed with the simplest
  choice,  
  ${p_i^{\mathsc h} = 0}$
  and
  $P_i = 0$.
Finally it is noted that the high-order contribution \eqref{eq:S:ansatz} 
  includes the viscous divergence contribution $\V F_{d3}$.
While this term vanishes for the exact solution,  
  it improves stability and accuracy with non-solenoidal approximations.


\section{Spatial discretization}
\label{sec:spatial discretization}

The semi-discrete SDC formulation developed in the previous section is 
  discretized in space using the discontinuous Galerkin spectral element
  method (DG-SEM) with nodal base functions \cite{Hesthaven2008a_SE}. 
Note that the following description is constrained to cuboidal domains.
This restriction serves only for convenience and can be lifted easily
  without affecting the proposed SDC method.

The description is organized as follows:
First, the necessary notation is introduced in 
  Sec.~\ref{sec:sdc:flow:dg:preliminaries}.
Section~\ref{sec:sdc:flow:dg:building blocks} defines the building blocks
  for designing the fully discrete SDC method, in particular
  the DG gradient and divergence functionals,
  the contributions to the time derivative,
  the interior penalty formulations of the pressure laplacian and the
  variable diffusion operator and the divergence/mass-flux stabilization.
Using these ingredients the discrete SDC predictor and corrector are
  composed in Sec.~\ref{sec:sdc:flow:dg:formulation}.
Following this,
Sec.~\ref{sec:sdc:flow:dg:solver} provides a short discussion on the 
  numerical evaluation of the discrete operators and, finally, 
  Sec.~\ref{sec:sdc:flow:dg:solver} summarizes the solution techniques
  and the implementation of the method.


\subsection{Preliminaries}
\label{sec:sdc:flow:dg:preliminaries}


First, the computational domain $\Omega$ is decomposed into
  $\Ne$ rectangular hexahedral elements to obtain the discrete domain
  \begin{equation}
    \label{eq:discrete domain}
    \Omega_h = \bigcup_{e=1}^{\Ne} \Omega^e
    \,.
  \end{equation}
Let 
  $\Gamma_h^{\mathsc i}$ denote the set of all interior 
    (including periodic) faces in $\Omega_h$
  and
  $\Gamma_h^{\d}$ the set of boundary faces.
The union of these sets defines the skeleton $\Gamma_h$.
For any interior face ${\Gamma^f \in \Gamma_h^{\mathsc i}}$
  there exist two adjoining elements $\Omega^-$ and $\Omega^+$ 
  with unit normal vectors $\V n^-$ and $\V n^+$, respectively.
The standard average and jump operators for element-wise continuous 
  functions $\V\phi$ of any dimension are defined as
  \begin{alignat}{2}
    \label{eq:average}
    & \avg{\V\phi}^f && = \tfrac{1}{2}(\V\phi^- + \V\phi^+)
    \,,
    \\
    \label{eq:jump:outer-product}
    & \jmpo{\V \phi}^f   && = \V n^- \V \phi^- + \V n^+ \V \phi^+
    \,,
  \end{alignat}
  where $\V\phi^\pm$ are the traces of the function from within $\Omega^\pm$.
A further jump operator, 
  involving the inner product with the normal vectors, 
  is introduced for vector or higher rank tensor functions:
  \begin{equation}
    \label{eq:jump:inner-product}
    \jmpi{\V\phi}^f = \V n^- \cdot \V\phi^- + \V n^+ \cdot \V\phi^+
    \,. 
  \end{equation}
These definitions are extended to boundary faces by assuming
  ${\V n^- = \V n = -\V n^+}$
  and providing exterior values 
  ${\V\phi^+ = \V\phi^{\d}}$
  on $\Gamma_h^{\d}$ 
  depending on boundary conditions
  \cite{Hesthaven2008a_SE}.
Hereafter, the index $f$ is dropped to indicate a quantity 
  that is defined on any face or a set of faces.

Let $\mathbb Q_{P}(\Omega^e)$ denote 
  the tensor-product space of all polynomials on $\Omega^e$
  with degree less or equal $P$ in each direction.
Glueing all element spaces together yields the global space 
  of element-wise polynomial, discontinuous functions 
  \begin{equation}
    \mathbb Q_{P} = \bigoplus_{\Omega^e \in \Omega_h} \mathbb Q_{P}(\Omega^e)
    \,.
  \end{equation}
This allows to define the ansatz spaces 
  \begin{gather}
    \mathbb V = [\mathbb Q_{P_v}]^3
    \quad 
    \text{and}
    \quad 
    \mathbb P = \mathbb Q_{P_p}
  \end{gather}
  for velocity and pressure, respectively.
Note that ${P_v > P_p}$ is required for inf-sup stability, see e.g. 
  \cite{Boffi2013a_SE, 
        John2016a_SE}. 
In the following the degree is set to 
  ${P_v = P}$   for velocity and
  ${P_p = P-1}$ for pressure.


\subsection{DG-SEM building blocks}
\label{sec:sdc:flow:dg:building blocks}


\subsubsection{Gradient and divergence functionals}
\label{sec:sdc:flow:dg:grad-div}

The semi-discrete SDC equations 
  (\ref{eq:predictor:1}--\ref{eq:predictor:4},
   \ref{eq:corrector:ansatz:1}--\ref{eq:S:ansatz},
   \ref{eq:H:ansatz})
  and, particularly, the time derivative contributions
  (\ref{eq:F_c}--\ref{eq:F_p})
  are composed of numerous gradient and divergence terms.
In the DG formulation, most of these terms can be expressed
  by generic functionals that are introduced in the following.
Consider the test functions
  ${q_h \in \mathbb P}$ 
  and 
  ${\V w_h \in \mathbb V}$.
The gradient functional of a scalar ${p_h \in \mathbb P}$ is given by
  \begin{equation}
    \label{eq:scalar-gradient-functional}
    \mathcal G_h(p_h, \V w_h)
       = - \int_{\Omega_h} (\nabla \cdot \V w_h) p_h \,\D \Omega
         + \int_{\Gamma_h} \jmpi{\V w_h} \avg{p_h} \D \Gamma
    \,.
  \end{equation}
Further,
  \begin{alignat}{2}
    \label{eq:vector-divergence-functional}
    & \mathcal D_h(\V v_h, q_h)
    &&  = - \int_{\Omega_h} \nabla q_h \cdot \V v_h \,\D \Omega
          + \int_{\Gamma_h} \jmpo{q_h} \cdot \avg{\V v_h} \D \Gamma
  \,,
  \\
    \label{eq:tensor-divergence-functional}
    & \mathcal D_h(\T\sigma_h, \V w_h)
    &&  = - \int_{\Omega_h} (\nabla \V w_h) : \T\sigma_h \,\D \Omega
          + \int_{\Gamma_h} \jmpo{\V w_h} : \avg{\T\sigma_h} \D \Gamma
  \end{alignat}
  define the divergence functionals for any
    vector 
      ${\V v_h \in \mathbb V}$ 
    and 
    second rank tensor
      ${\T\sigma_h \in \mathbb V \otimes \mathbb V}$,
    respectively.
Boundary conditions are considered by providing proper exterior values, 
  as will be detailed below.
Based on these functionals 
  the discrete 
  gradient ${\nabla_h p_h}$ and
  divergence ${\nabla_h \cdot \V v_h}$
  are introduced such that
  \begin{alignat}{3}
    & \int_{\Omega_h} \V w_h \cdot \nabla_h p_h \,\D \Omega
    && = \mathcal G_h(p_h, \V w_h)
    &\quad& \forall q_h \in \mathbb P 
    \,,
    \\
    & \int_{\Omega_h} q_h \nabla_h \cdot \V v_h \,\D \Omega
    && = \mathcal D_h(\V v_h, q_h)
    &\quad& \forall \V w_h \in \mathbb V
    \,.
  \end{alignat}
These forms were already advocated by \citet{Krank2017a_TI}, 
  who found that the underlying partial integration 
  improves the robustness of their projection method,
  when combined with central fluxes (averages) of $\V v_h$ and $p_h$ 
  across the element boundaries.

\subsubsection{Time derivative}
\label{sec:sdc:flow:dg:time derivative}

The time derivative comprises the discrete counterparts of the
  convection, diffusion, pressure and forcing terms introduced in
  (\ref{eq:F}\,--\,\ref{eq:F_d}).
For the convection term $\V F_{h,c}$ application of the local Lax-Friedrichs 
  flux leads to
  \begin{equation}
    \label{eq:Fhc}
    \begin{aligned}
      \int_{\Omega_h} \V w_h \cdot \V F_{h,c}\,\D \Omega
         ~=~ & - \int_{\Omega_h} \nabla\V w_h : \V v_h\V v_h\,\D \Omega
      \\     & + \int_{\Gamma_h} \jmpo{\V w_h} : 
                               \big( \avg{\V v_h\V v_h}
                                   + \hat v_n \jmpo{\V v_h}
                               \big) \,\D \Gamma
      \quad \forall \V w_h \in \mathbb V
      \,,
    \end{aligned}
  \end{equation}
  where
  ${\hat v_n = \max\big( |\V n \cdot \V v_h^-|,|\V n \cdot \V v_h^+| \big)}$.
  It is worth noting that this flux adds artificial dissipation to the 
  method, which disappears, however, if the velocity jumps $\jmpo{\V v_h}$
  vanish, see \cite{Hesthaven2008a_SE}.

The diffusive and pressure terms are based on the divergence and gradient
  functionals, i.e.
  \begin{alignat}{2}
    \label{eq:Fhd1}
    & \int_{\Omega_h} \V w_h \cdot \V F_{h,d1}\,\D \Omega
    && ~=~ \mathcal D_h(\nu_h \nabla \V v_h, \V w_h)
    \,,
    \\
    \label{eq:Fhd2}
    & \int_{\Omega_h} \V w_h \cdot \V F_{h,d2}\,\D \Omega
    && ~=~ \mathcal D_h(\nu_h \transpose{(\nabla \V v_h)}, \V w_h)
    \,,
    \\
    \label{eq:Fhd3}
    & \int_{\Omega_h} \V w_h \cdot \V F_{h,d3}\,\D \Omega
    && ~=~ -\chi\, \mathcal G_h(\nu_h \nabla \cdot \V v_h, \V w_h)
    \,,
    \\
    \label{eq:Fhp}
    & \int_{\Omega_h} \V w_h \cdot \V F_{h,p}\,\D \Omega
    && ~=~ \mathcal G_h(p_h, \V w_h)
  \end{alignat}
  for all ${\V w_h \in \mathbb V}$.
Similarly, the forcing term follows from
  \begin{equation}
    \label{eq:fh}
    \int_{\Omega_h} \V w_h \cdot \V f_{h}\,\D \Omega
      = \int_{\Omega_h} \V w_h \cdot \V f\,\D \Omega
    \quad \forall \V w_h \in \mathbb V
    \,,
  \end{equation}
  which is equivalent to element-wise $L^2$ projection.

\subsubsection{Laplace and viscous diffusion operators}
\label{sec:sdc:flow:dg:diffusion}

The Laplacian occurring in the pressure equations such as
  \eqref{eq:pressure problem}
  and the viscous diffusion operator in
  \eqref{eq:predictor:3} and \eqref{eq:corrector:ansatz:3}
  are discretized using the symmetric interior penalty (SIP) 
  method \cite{Arnold2001a_SE}. 

Application to the pressure Laplacian $\nabla^2 p$ yields
  \begin{equation}
    \label{eq:L_h}
    \begin{aligned}
      \mathcal L_h(p_h, q_h)
      = & - \int_{\Omega_h} \nabla q_h \cdot \nabla p_h \,\D \Omega
      \\
        & + \int_{\Gamma_h} (
                              \jmpo{q_h} \cdot \avg{\nabla p_h}
                            + \avg{\nabla q_h} \cdot \jmpo{p_h}
                            ) \,\D \Gamma
      \\    
        & - \int_{\Gamma_h} \mu_p\jmpo{q_h} \cdot \jmpo{p_h} \,\D \Gamma
      \,,
    \end{aligned}
  \end{equation}
  for ${p_h,q_h \in \mathbb P}$.
  The penalty parameter is defined as ${\mu_p = \mu(P_p)}$ with
  \begin{equation}
    \mu(P) = \mu_{\star} \frac{P(P+1)}{2} \avg{\frac{1}{\Delta x_n}},
  \end{equation}
  where 
  $\Delta x_n$ is the mesh spacing normal to the face
  and
  ${\mu_{\star} > 1}$
  a constant parameter 
  \cite{Stiller2016b_SE}.
On Neumann boundaries the face averages and jumps are given by
  \begin{equation}
    \V n \cdot \avg{\nabla p} = (\d_n p)_b,
    \quad
    \jmpo{p_h} = 0
    \,.
  \end{equation}
Note that the latter implies ${p^+ = p^-}$.

For the diffusion operator 
  ${\lambda \V v - \nabla\cdot\nu\nabla\V v}$
  with constant $\lambda$
  and variable $\nu$
  the SIP method gives
  \begin{equation}
    \label{eq:V_h}
    \begin{aligned}
      \mathcal V_h(\lambda, \nu_h; \V v_h, \V w_h)
      & = \int_{\Omega_h} \lambda \V w_h \cdot \V v_h
        + \int_{\Omega_h} \nu_h \nabla \V w_h : \nabla \V v_h \,\D \Omega
      \\
      & - \int_{\Gamma_h} (
                              \jmpo{\V w_h} : \avg{\nu_h\nabla \V v_h}
                            + \avg{\nu_h\nabla \V w_h} : \jmpo{\V v_h}
                           ) \,\D \Gamma
      \\
      & + \int_{\Gamma_h} 
            \mu_v\,\hat\nu\,\jmpo{\V w_h} : \jmpo{\V v_h} \,\D \Gamma
    \end{aligned}
  \end{equation}
  where 
  ${\mu_v = \mu_{\star}(P_v)}$
  and
  ${\hat \nu = \max(\nu_h^-, \nu_h^+)}$.
Dirichlet boundary conditions are weakly imposed by setting
  \begin{equation}
    \avg{\V v} = \V v_b,
    \quad
    \avg{\nu_h \nabla\V v_h} = (\nu_h \nabla\V v_h)^-
    \,.
  \end{equation}
The first of these relations is equivalent to
  ${\V v^+ = 2 \V v_b - \V v^-}$
  and thus also defines the jump
  ${\jmpo{\V v_h}}$.
Remarkably, $\mathcal V_h$ does not couple across $\V v_h$ 
  such that the corresponding diffusion problems can be solved 
  component by component as long as $\nu_h$ and the RHS do not
  depend on the solution.

\subsubsection{Divergence/mass-flux stabilization.}
\label{sec:sdc:flow:dg:stabilization}

The discrete projection steps are augmented with the penalty functional
  \begin{equation}
    \label{eq:J_h}
    \mathcal J_h(\V v_h, \V w_h)
      = \int_{\Omega_h} \tau_d (\nabla\cdot\V w_h)(\nabla\cdot\V v_h)\,\D \Omega
      + \int_{\Gamma_h} \tau_j \jmpi{\V w_h} \jmpi{\V v_h} \,\D \Gamma
    \,.
  \end{equation}
This functional was introduced by \citet{Joshi2016a_SE} in frame of a
  post-processing technique for stabilizing pressure-correction methods 
  for incompressible inviscid flow.
It has no counterpart in the differential formulation, but vanishes for
  continuous, element-wise divergence-free vector fields $\V v_h$.
In the general case, the first part of $\mathcal J_h$ penalizes the
  divergence of $\V v_h$ within elements and the second part 
  jumps of the normal flux across faces.
\citet{Akbas2018a_SE} recently proved that both parts are required
  for pressure robustness.

The divergence penalty functional \eqref{eq:J_h} has been applied with
  projection methods as well as coupled methods
  \cite{Joshi2016a_SE,Krank2017a_TI,Akbas2018a_SE,Fehn2019a_SE}.
In these studies, various expressions have been proposed for 
  the stabilization parameters $\tau_d$ and $\tau_j$.
The present work follows \cite{Akbas2018a_SE} by setting
  \begin{gather}
    \tau_d = \tau_\star \,\nu_{\mathrm{ref}}
    \qquad
    \text{and}
    \qquad
    \tau_j = \frac{\tau_\star}{\Delta x_n} \,\nu_{\mathrm{ref}}
    \,,
  \end{gather}
  where 
    $\tau_\star$ is a positive constant, 
    $\Delta x_n$ the mesh spacing in the normal direction 
  and
    $\nu_{\mathrm{ref}}$ a reference value of viscosity.

\subsection{DG-SEM formulation of the SDC method}
\label{sec:sdc:flow:dg:formulation}

\subsubsection{Predictor}
\label{sec:sdc:flow:dg:predictor}

Application of the above building blocks leads to discrete versions
  of the predictor and the corrector.
Both resemble their semi-discrete precursors, but are further transformed
  here to reflect the actual course of computation.

Starting from $\V v^0_{h,0}$ at time $t_0$ the predictor sweeps
  through all subintervals $i$, performing the following steps:
\begin{enumerate}

\item
  \emph{Extrapolation:}
  \begin{equation}
    \label{eq:predictor:dg:extrapolation}
    \V v'^{\,0}_{h,i} 
         = \V v^0_{h,i-1}
           + \Delta t_i 
             \big[ \V F_{h,c}(\V v^0_{h,i-1}) 
                 + \V F_{h,d}(\nu^0_{h,i-1}, \V v^0_{h,i-1})
                 + \V f_{h,i}
             \big]
             \,.
  \end{equation}

\item
  \emph{First projection:} 
  Determine pressure ${p''^{\,0}_{h,i} \in \mathbb P}$ by solving
  \begin{equation}
    \label{eq:predictor:dg:projection1:pressure}
    \mathcal L_h(p''^{\,0}_{h,i}, q_h) 
      = -\frac{1}{\Delta t_i}\,\mathcal D_h(\V v'^{\,0}_{h,i}, q_h)
      \quad
      \forall q_h \in \mathbb P
  \end{equation}
  with Neumann boundary conditions
  \begin{equation}
    \V n \cdot \nabla p''^{\,0}_{h,i}|_{\d\Omega_h} = 
      \frac{1}{\Delta t_i}\,
      \V n \cdot (\V v'^{\,0}_{h,i}|_{\d\Omega_h} - \V v_{\textrm b}(t_i))
      \,.
  \end{equation}
  Subsequently compute the intermediate velocity 
  ${\V v''^{\,0}_{h,i} \in \mathbb V}$
  such that
  \begin{equation}
    \label{eq:predictor:dg:projection1:correction}
    \int_{\Omega_h} 
       \V w_h \cdot \frac{\V v''^{\,0}_{h,i} - \V v'^{\,0}_{h,i}}{\Delta t_i}
       \,\D \Omega
     + \mathcal J_h(v''^{\,0}_{h,i}, \V w_h)
     + \mathcal G_h(p''^{\,0}_{h,i}, \V w_h)
     = 0
     \quad
     \forall \V w_h \in \mathbb V
  \end{equation}
  with homogeneous Neumann conditions (constant extrapolation)
  on $\d\Omega_h$.

\medskip
\item
  \emph{Diffusion:}
  Find ${\V v'''^{\,0}_{h,i} \in \mathbb V}$ such that for all 
  ${\V w_h \in \mathbb V}$
  \begin{equation}
    \label{eq:predictor:dg:diffusion}
    \begin{aligned}
      \mathcal V_h(\Delta t_i^{-1}, \nu^0_{h,i-1}; \V v'''^{\,0}_{h,i}, \V w_h)
        & = \int_{\Omega_h} 
              \V w_h \cdot \frac{\V v''^{\,0}_{h,i}}{\Delta t_i}
              \,\D \Omega 
      \\
        & - \int_{\Omega_h} 
              \V w_h \cdot (\V F_{h,d1} + c_{\chi}\V F_{h,d3})_{i-1}^0
              \,\D \Omega
     \end{aligned}
  \end{equation}
  with Dirichlet conditions
  \begin{equation}
    \V v'''^{\,0}_{h,i}|_{\d\Omega_h} = \V v_{\textrm b}(t_i)
    \,.
  \end{equation}

\item
  \emph{Final projection:} 
  Solve
  \begin{equation}
    \label{eq:predictor:dg:projection2:pressure}
    \mathcal L_h(p^0_{h,i} - p''^{\,0}_{h,i}, q_h) 
       = -\frac{1}{\Delta t_i}\,
         \mathcal D_h(\V v'''^{\,0}_{h,i}, q_h)
    \quad
    \forall q_h \in \mathbb P
  \end{equation}
  for ${p^0_{h,i} \in \mathbb P}$ 
  and, subsequently, 
  \begin{equation}
    \label{eq:predictor:dg:projection2:correction}
    \int_{\Omega_h} 
       \V w_h \cdot \frac{\V v^0_{h,i} - \V v'''^{\,0}_{h,i}}{\Delta t_i}
       \,\D \Omega
     + \mathcal J_h(\V v^0_{h,i}, \V w_h)
     + \mathcal G_h(p^0_{h,i} - p''^{\,0}_{h,i}, \V w_h)
     = 0
     \quad
     \forall \V w_h \in \mathbb V
  \end{equation}
  to obtain the velocity ${\V v^0_{h,i} \in \mathbb V}$.
  For both problems,
  \eqref{eq:predictor:dg:projection2:pressure} 
  as well as
  \eqref{eq:predictor:dg:projection2:correction}, 
  homogeneous Neumann conditions are imposed.
  
\end{enumerate}

\subsubsection{Corrector}
\label{sec:sdc:flow:dg:corrector}

The DG-SEM formulation of the corrector resembles that of the predictor.
It is summarized below, skipping specifications of function spaces and
  homogeneous boundary conditions for brevity.

For 
  ${k=0,\dots,K-1}$
  sweep
  through all subintervals $i$ and perform the following steps:
\begin{enumerate}[start=0]

\item
  \emph{Low- and high-order contributions:}
  \begin{alignat}{2}
  & \V H^{k}_{h,i} && = 
      \Delta t_i 
        \big[ \V F_{h,c,i-1} 
            + \V F_{h,d1}(\nu_{h,i-1}, \V v'''_{h,i})
            + \V F_{h,d2,i-1}
            + c_{\chi}\V F_{h,d3,i-1}
        \big]^k
        \,,
  \\
  & \V S^{k}_{h,i} && =
      \Delta t \sum_{j=0}^M w^{\mathsc s}_{i,j} 
        \big[ \V F_{h,c,i}^k
            + \V F_{h,d,i}^k
            + \V f_{h,i}
        \big]
        \,.
  \end{alignat}

\item
  \emph{Extrapolation:}
  \begin{equation}
    \label{eq:corrector:dg:extrapolation}
    \begin{aligned}
      \V v'^{\,k+1}_{h,i} 
         & = \V v^{k+1}_{h,i-1}
           + \Delta t_i 
               \big[ \V F_{h,c,i-1}^{k+1} 
                   + \V F_{h,d1}(\nu_{h,i-1}^{k+1}, \V v^{k}_{h,i})
                   + \V F_{h,d2,i-1}^{k+1}
                   + \V F_{h,d3,i-1}^{k+1}
               \big]   
      \\ & - \V H^{k}_{h,i}
           + \V S^{k}_{h,i}
      \,.
    \end{aligned}
  \end{equation}

\item
  \emph{First projection:} 
  \begin{gather}
    \label{eq:corrector:dg:projection1:pressure}
    \mathcal L_h(\tilde p''^{\,k+1}_{h,i}, q_h) =
    - \frac{1}{\Delta t_i}\,\mathcal D_h(\V v'^{\,k+1}_{h,i}, q_h)
    \,,
  \\
    \V n \cdot \nabla \tilde p''^{\,k+1}_{h,i}|_{\d\Omega_h} = 
      \frac{1}{\Delta t_i}\,
      \V n \cdot (\V v'^{\,k+1}_{h,i}|_{\d\Omega_h} - \V v_{\textrm b}(t_i))
      \,,
  \\
    \label{eq:corrector:dg:projection1:correction}
    \int_{\Omega_h} 
        \V w_h \cdot \frac{\V v''^{\,k+1}_{h,i} - \V v'^{\,k+1}_{h,i}}{\Delta t_i}
        \,\D \Omega
    + \mathcal J_h(v''^{\,k+1}_{h,i}, \V w_h)
    + \mathcal G_h(\tilde p''^{\,k+1}_{h,i}, \V w_h)
    = 0
    \,.
  \end{gather}

\item
  \emph{Diffusion:}
  \begin{gather}
    \label{eq:corrector:dg:diffusion}
    \begin{aligned}
      \mathcal V_h(\Delta t_i^{-1}, \nu^{k+1}_{h,i-1}; \V v'''^{\,k+1}_{h,i}, \V w_h)
        & = \int_{\Omega_h} 
              \V w_h \cdot \frac{\V v''^{\,k+1}_{h,i}}{\Delta t_i}
              \,\D \Omega 
     \\
        & - \int_{\Omega_h} 
              \V w_h \cdot \big[ \V F_{h,d1}(\nu^{k+1}_{h,i-1}, \V v^{k}_{h,i})
                           + c_{\chi}\V F_{h,d3,i-1}^{k+1}
                           \big]
              \,\D \Omega
    \end{aligned}
  \end{gather}
  with Dirichlet conditions \,
  \(
    \V v'''^{\,k+1}_{h,i}|_{\d\Omega_h} = \V v_{\textrm b}(t_i)
    \,.
  \)

\medskip
\item
  \emph{Final projection:} 
  \begin{gather}
    \label{eq:corrector:dg:projection2:pressure}
    \mathcal L_h(\tilde p^{k+1}_{h,i} - \tilde p''^{\,k+1}_{h,i}, q_h) 
    = -\frac{1}{\Delta t_i}\,
       \mathcal D_h(\V v'''^{\,k+1}_{h,i}, q_h)
    \,,
  \\
    \label{eq:corrector:dg:projection2:correction}
    \int_{\Omega_h} 
      \V w_h \cdot \frac{\V v^{k+1}_{h,i} - \V v'''^{\,k+1}_{h,i}}{\Delta t_i}
      \,\D \Omega
    + \mathcal J_h(\V v^{k+1}_{h,i}, \V w_h)
    + \mathcal G_h(\tilde p^{k+1}_{h,i} - \tilde p''^{\,k+1}_{h,i}, \V w_h)
    = 0
    \,.
  \end{gather}
  
\end{enumerate}

\subsubsection{Base functions and numerical quadrature}
\label{sec:sdc:flow:dg:quadrature}

The discrete solution is approximated by means of 
  tensor-product Lagrange bases constructed from GLL points
  of degree
  $P_p$ for the pressure and
  $P_v$ for the velocity and all remaining variables.
Integrals are evaluated numerically 
  with GLL quadrature on the collocation points 
  except for
  the convection term \eqref{eq:Fhc} and
  the functional \eqref{eq:scalar-gradient-functional},
  which are integrated using
  ${\big\lceil\frac{3P_v}{2}\big\rceil\!+\!1}$ 
  and
  ${P_v\!+\!1}$ 
  points, respectively.
This choice avoids aliasing errors and preserves the equivalence 
  between the gradient and divergence functionals, i.e.
  \begin{equation}
    \mathcal G_h(p_h, \V v_h) = -\mathcal D_h(\V v_h, p_h)
    \,.
  \end{equation}
According to \citet{Maday1990a_SE} 
  optimal convergence with variable viscosity requires
  ${P_v\!+\!2}$ 
  Lobatto points for integrating the diffusion term \eqref{eq:V_h},
  as opposed to ${P_v\!+\!1}$ which are actually used.
However, elevating the quadrature order would also increase the cost 
  of solving the diffusion problems and is therefore postponed 
  to future studies.


\subsection{Solution methods and implementation}
\label{sec:sdc:flow:dg:solver}


The pressure equations 
  (\ref{eq:predictor:dg:projection1:pressure},
   \ref{eq:predictor:dg:projection2:pressure},
   \ref{eq:corrector:dg:projection1:pressure})
  and
  (\ref{eq:corrector:dg:projection2:pressure})
  are solved by means of a Krylov-accelerated polynomial multigrid technique
  using an element-based overlapping Schwarz method for smoothing
  \cite{Stiller2016b_SE,Stiller2017a_SE}.
To cope with variable coefficients in diffusion problems 
  (\ref{eq:predictor:dg:diffusion},
   \ref{eq:corrector:dg:diffusion})
  the Schwarz smoother was extended
  by adopting the linearization strategy developed in 
  \cite{Stiller2016a_SE}
  for continuous spectral elements.
The projection steps
  (\ref{eq:predictor:dg:projection1:correction},
   \ref{eq:predictor:dg:projection2:correction},
   \ref{eq:corrector:dg:projection1:correction})
  and
  (\ref{eq:corrector:dg:projection2:correction})
  are solved with a diagonally preconditioned conjugate gradient method
  \cite{Shewchuk1994a_KR}.
The SDC method and all examples presented in this paper are 
  implemented in the high-order spectral-element techniques library
  \emph{HiSPEET}\footnote{%
  \emph{HiSPEET} is freely available as a git repository at
  \url{fusionforge.zih.tu-dresden.de/projects/hispeet}.
  Since the library is still in an early stage of development, 
  the reader is encouraged to contact the author for further
  instructions.}
All solver components are parallelized with MPI and 
  exploit SIMD techniques for accelerating the element operators on CPUs 
  \cite{Huismann2020a_SE}.


\section{Numerical experiments}
\label{sec:numerical experiments}

In the following, numerical results are presented for various test cases
  including those with constant viscosity as well as several scenarios 
  with variable viscosity.
For the pressure and diffusion problems the multigrid solver is used with
  a relative tolerance of 10\textsuperscript{\textminus 12} and 
  an absolute tolerance of 10\textsuperscript{\textminus 14}.
The preconditioned conjugate gradient method in the projection steps 
  is terminated after 10 iterations or reaching 
  a relative tolerance of 10\textsuperscript{\textminus 10}.
Unless stated otherwise, the parameter of the divergence/mass-flux 
  stabilization is set to ${\tau_\star=1}$.
The velocity error $\varepsilon_v$ is computed as the RMS value over all 
  mesh points at the end time $T$.
Correspondingly, 
  $\varepsilon_p$ denotes the pressure error and
  $\varepsilon_{\mathrm{div}}$ the divergence error,
  i.e., the RMS value of ${\nabla_h \cdot \V v_h}$.


\subsection{Traveling Taylor-Green vortex with constant viscosity}

The first test problem is adopted from \citet{Minion2018a_TI}.
It is a Taylor-Green vortex traveling through the two-dimensional domain
  ${\Omega^{\mathrm{2D}} = [ -\sfrac{1}{2}, \sfrac{1}{2} ]^2}$.
The exact solution is given by
  \begin{alignat}{2}
    & v_x^{\mathrm{ex}}(x,y,t) 
    &&= 1 + \sin \big( 2\pi ( x                - t )\big) 
            \cos \big( 2\pi ( y - \tfrac{1}{8} - t )\big) 
            \exp( -8\pi^2\nu t ) 
    \,,
    \\
    & v_y^{\mathrm{ex}}(x,y,t) 
    &&= 1 - \cos \big( 2\pi ( x                - t )\big)
            \sin \big( 2\pi ( y - \tfrac{1}{8} - t )\big)
            \exp( -8\pi^2\nu t ) 
    \,,
    \\
    & p^{\mathrm{ex}}(x,y,t)   
    &&=     \tfrac{1}{4}\big[ \cos \big(4\pi (x                - t)\big)
                            + \cos \big(4\pi (y - \tfrac{1}{8} - t)\big) 
                            \big]
            \exp(-16\pi^2\nu t)
    \,.
  \end{alignat}
There are no external sources, i.e., ${\V f = 0}$.
Based on this problem \citet{Minion2018a_TI} defined two test cases:
  Example\,1 
    with 
    ${\nu = 0.02}$ and
    periodic conditions in both directions, 
  and
  Example\,2
    with
    ${\nu = 0.01}$,
    periodic conditions in $x$-direction and
    Dirichlet conditions in $y$-direction.
The final time is ${T=0.25}$ in both cases.
For the present study these cases are extended into three dimensions 
  by assuming ${v_z = 0}$ and periodicity in the $z$-direction.
Following \cite{Minion2018a_TI} the domain is discretized with
  8\texttimes8 elements of degree ${P=10}$ in the $x$-$y$ plane 
  and one element layer in the $z$-direction.
Time integration is performed with the SDC method using ${M=3}$ 
  subintervals.
The corresponding numerical studies are labeled 
  PPP for the periodic case and 
  PDP for the case with Dirichlet conditions in the $y$-direction.

\begin{figure}
  \includegraphics[height=67mm]{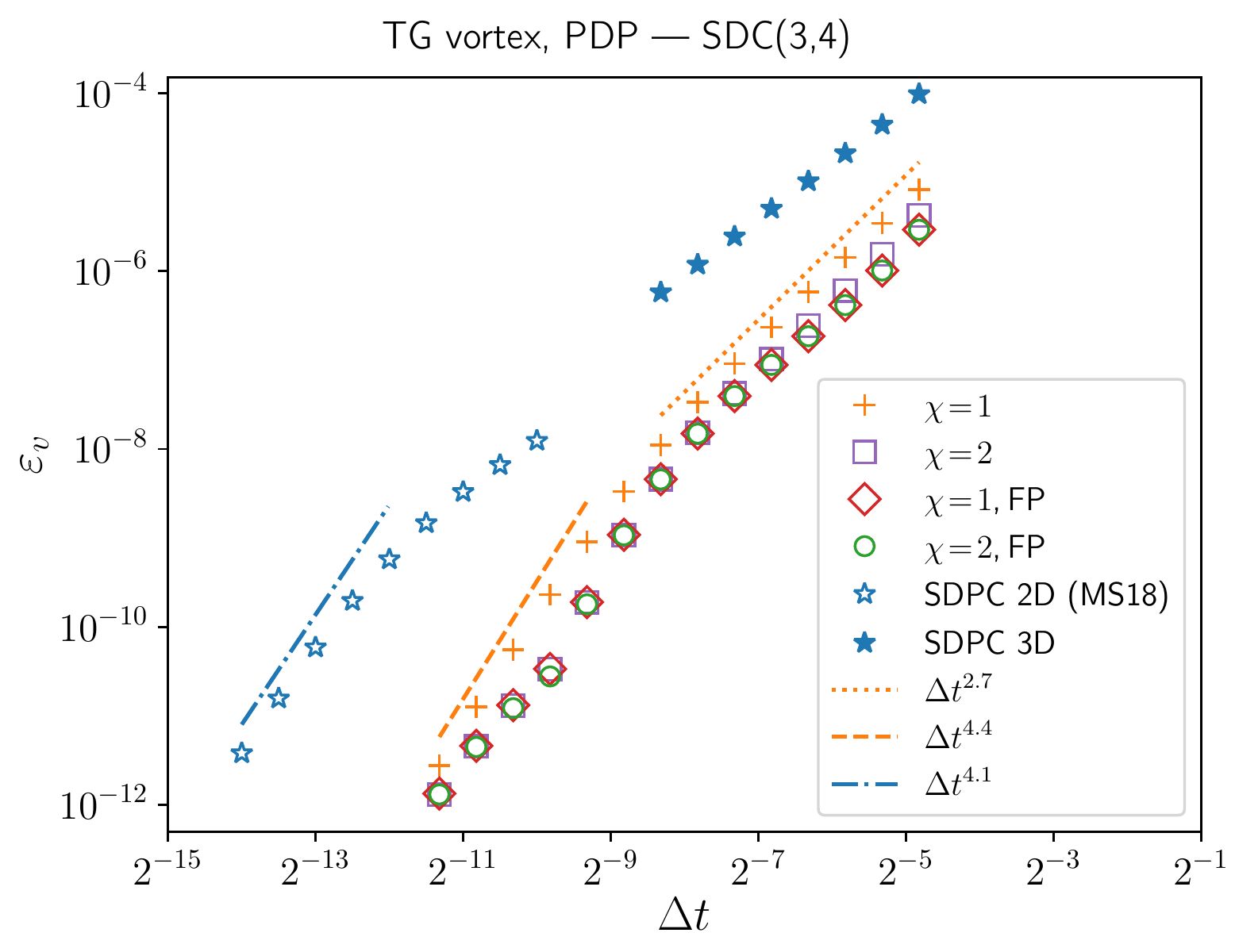}
  \caption{%
    Velocity error of SDC with 
    ${M\!=\!3}$ subintervals and ${K\!=\!4}$ correction sweeps
    for the Taylor-Green vortex with Dirichlet conditions in the $y$-direction.
    The choices
    ${\chi\!=\!1}$ 
    and
    ${\chi\!=\!2}$ 
    refer to the method based on standard and rotational
    velocity correction, respectively, and 
    FP indicates the application of the final projection step.
    SDPC 2D (MS18) represents the two-dimensional results of \citet{Minion2018a_TI} and
    SDPC 3D corresponding three-dimensional results obtained in the present study.
    The latter extends only down to $\Delta t \approx 2^{-8}$, as the
    method became increasingly unstable with smaller time steps.
    (SDPC 2D results courtesy of Michael Minion)
    \label{fig:tg_pdp_sdc_3-4_comparison}
    }
\end{figure}

Before considering the studies in detail, selected variants of the
  present SDC method are compared to each other and to the SDPC 
  method of \citet{Minion2018a_TI}.
Figure \ref{fig:tg_pdp_sdc_3-4_comparison}
  shows the velocity error after ${K\!=\!4}$ correction sweeps.
As the most striking feature, the original 2D results of
  \cite{Minion2018a_TI} exhibit an error up to 1000 times greater 
  than all variants of the present SDC method for a given 
  step size $\Delta t$.
Given this unexpected deviation, a 3D version of SDPC was established,
  albeit using an IMEX Euler corrector, and applied to the test case.
However, the 3D SDPC also failed to match the present SDC method
  and, unfortunately, suffered from instability with time steps 
  ${\Delta t < 2^{-9}}$.
Possibly, its stability can be improved by using a DIRK sweeps
  as proposed in \cite{Weiser2015a_TI}.
This is, however, out of the scope of the present work.
Among the investigated variants of the proposed SDC method, the
  one based on the standard velocity correction scheme (${\chi=1}$)
  with no final projection shows the most consistent behavior.
For time steps larger than $2^{-9}$ it attains a convergence rate of 
  only $2.7$, which is scarcely more than half of the expected order 
  of 5.
This order reduction is absent in the purely periodic case and, hence,
  attributed to imposing unsteady Dirichlet conditions.
With smaller steps the convergence rate increases to 4.4, which is
  only about half an order less than the optimum.
Executing the final projection (FP) improves the accuracy throughout, but
  leads to a less regular convergence pattern. 
The rotational scheme achieves a similar improvement even without
  the final projection.
Since these methods differ only in the splitting scheme used for
  approximating the IMEX Euler method, it can be concluded that 
  minimizing the splitting error is of crucial importance for the 
  overall stability and accuracy.

\begin{figure}
\subfloat[Periodic conditions in $y$-direction.]
  {\includegraphics[height=67mm]{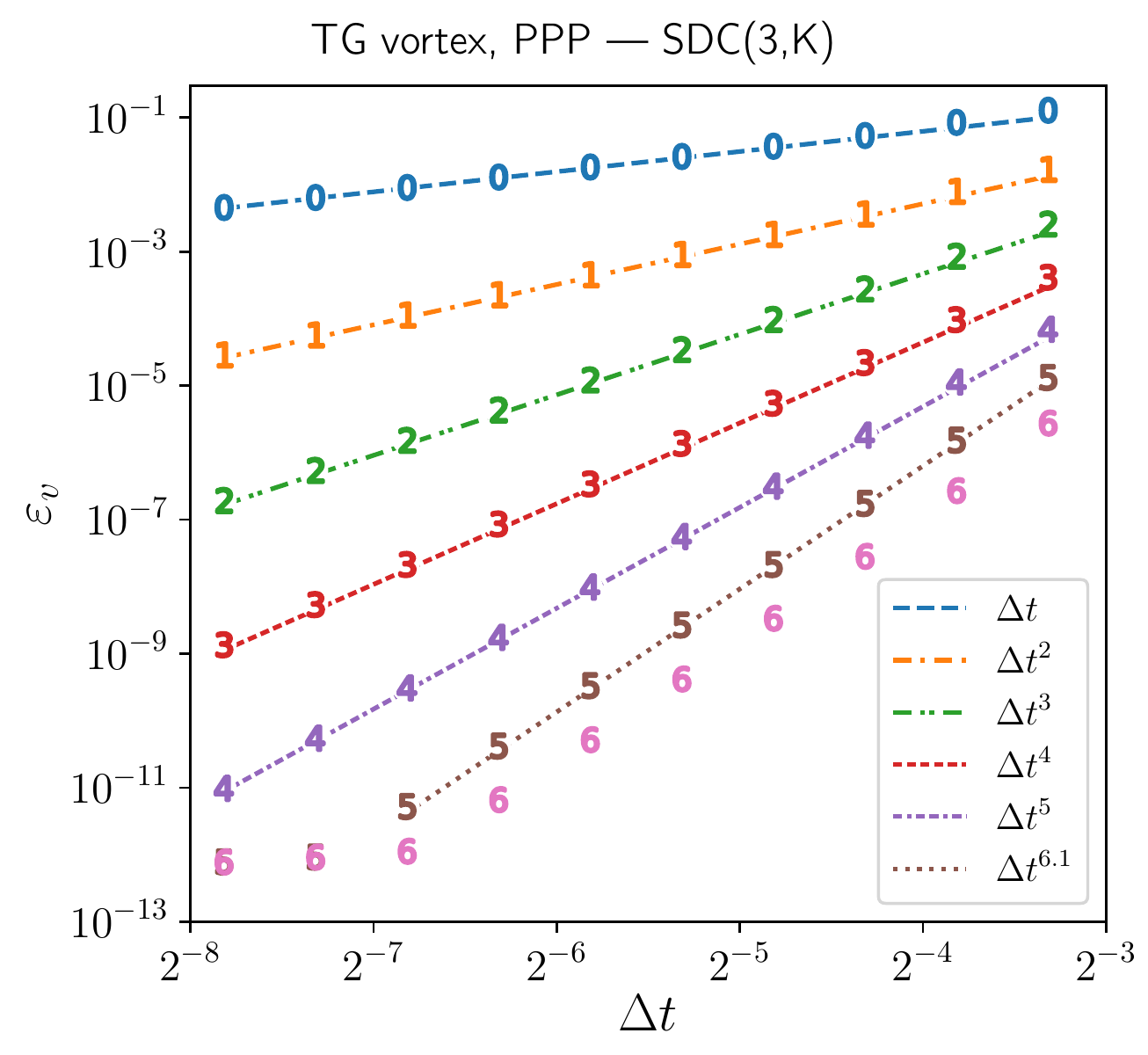}
   \label{fig:tg_ppp_sdc_3-v}}
\subfloat[Dirichlet conditions in $y$-direction.]
  {\includegraphics[height=67mm]{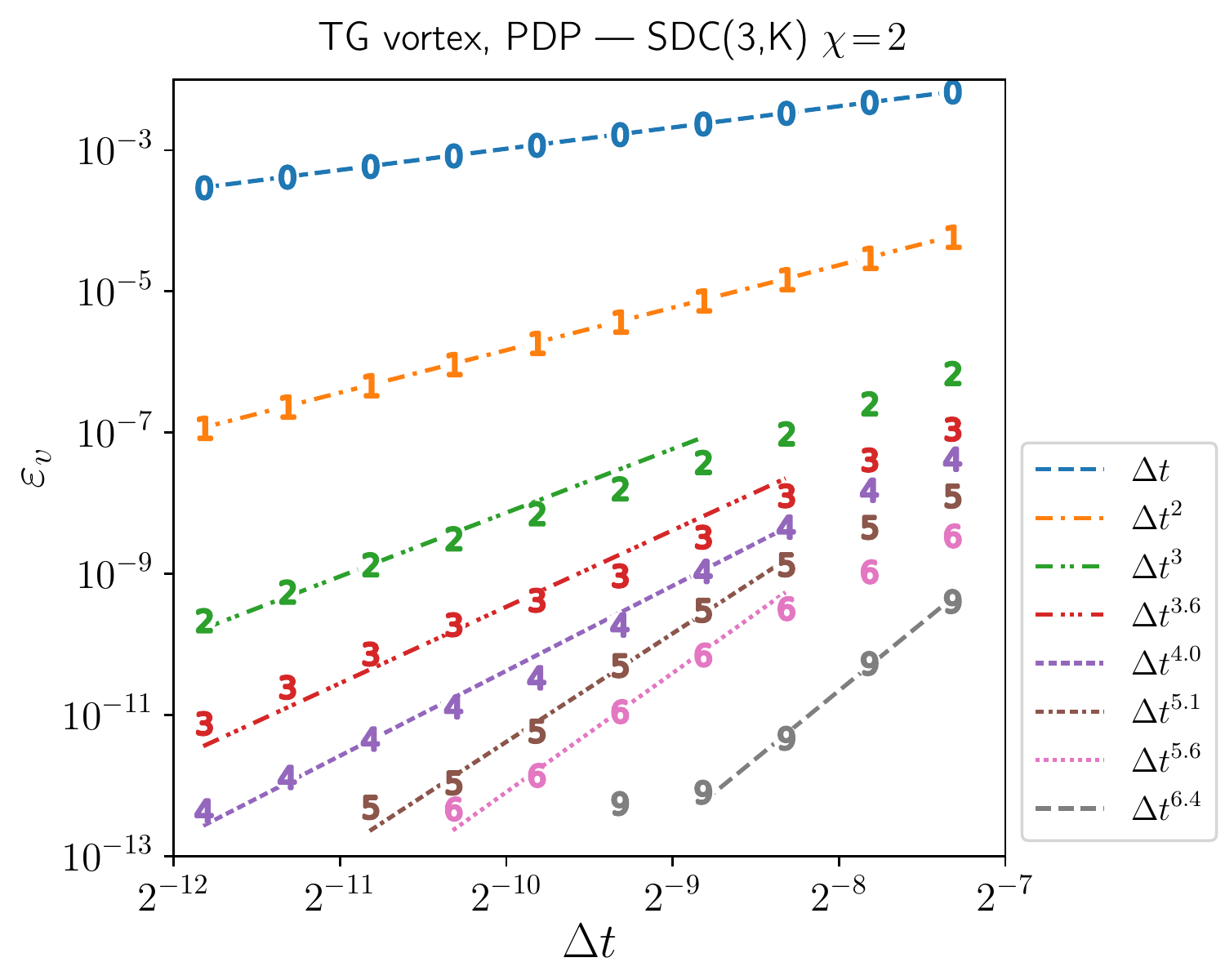}
   \label{fig:tg_pdp_sdc_3-v}}
\caption{%
  Velocity error for the Taylor-Green vortex with ${M=3}$ subintervals, 
  ${\chi\!=\!2}$ and no final projection.
  Labels indicate the number of correction sweeps and
  lines the approximate slope.
  \label{fig:tg_sdc_3-v}
  }
\end{figure}

\begin{figure}
  \includegraphics[height=67mm]{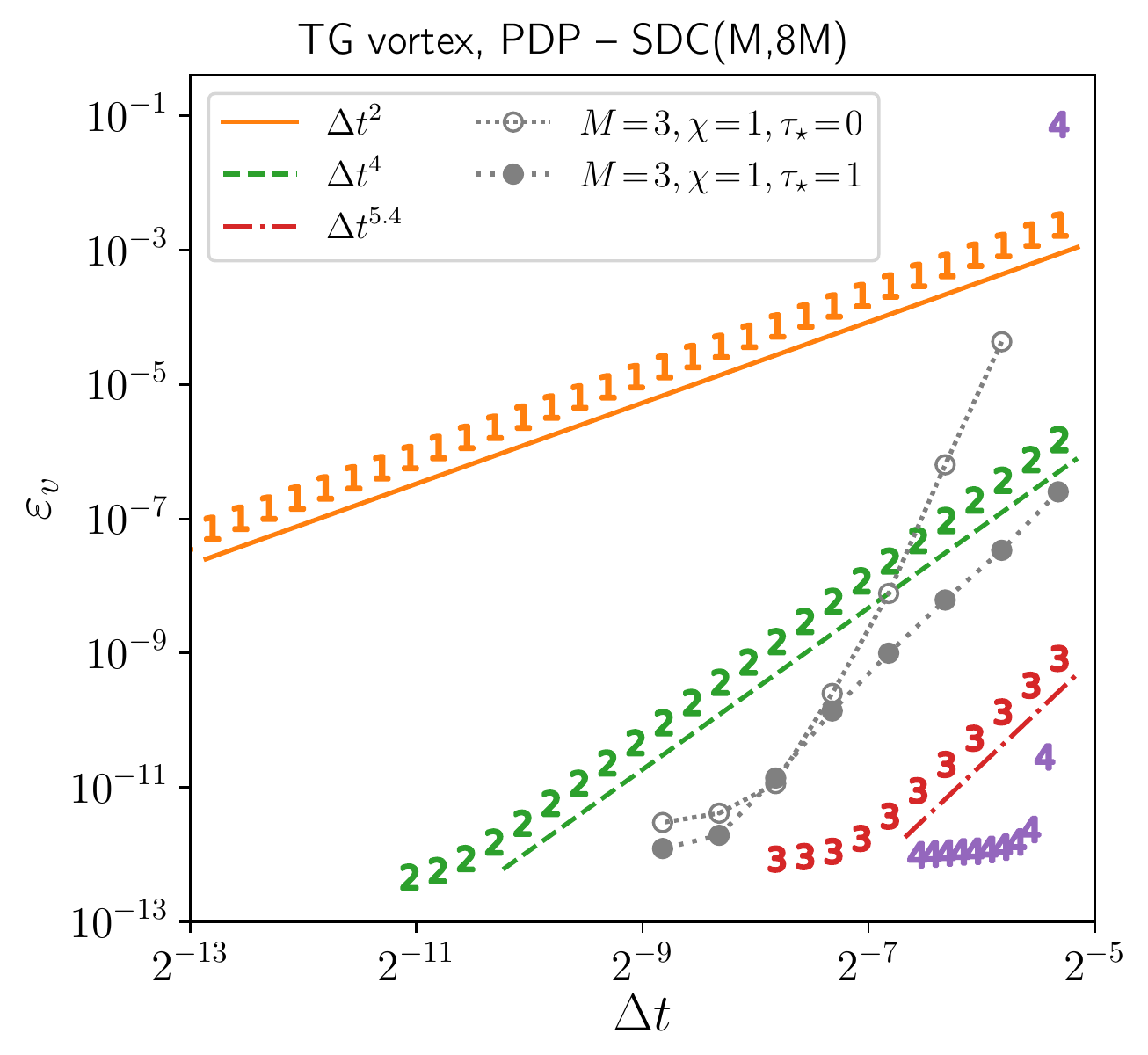}
\caption{%
  Velocity error of the converged SDC method for the Taylor-Green 
  vortex with Dirichlet conditions in the $y$-direction. 
  Labeled graphs refer to the rotational scheme (${\chi\!=\!2}$) 
  with divergence/mass-flux stabilization ($\tau_\star\!=\!1$) 
  and final projection. 
  Digits indicate the number $M$ of subintervals.
  The dashed lines correspond to the standard scheme (${\chi\!=\!1}$)
  with or without stabilization 
  ($\tau_\star\!=\!0$ or $\tau_\star\!=1$) 
  and no final projection.
  \label{fig:tg_sdc_m-v}
  }
\end{figure}

Figure \ref{fig:tg_sdc_3-v} presents the results for cases PPP and
  PDP obtained with the rotational method for a different number of
  correction sweeps, ranging from ${K=0}$ to $9$.
For PPP the convergence rate grows by one with each 
  correction until the maximum of $2M$ is reached
  (Fig.~\ref{fig:tg_ppp_sdc_3-v}).
Hence, the method shows the optimal convergence behavior
  in the periodic case.
This 
  can be explained by the lack of a splitting error
  and  
  was also observed in \cite{Minion2018a_TI}.
In case PDP
  the imposition of Dirichlet conditions 
  causes an order reduction which manifests
  in a larger error and a flatter slope for identical $K$
  in comparison to PPP
  (Fig.~\ref{fig:tg_pdp_sdc_3-v}).
Increasing the number of sweeps to ${K=9}$ yields a further error
  reduction and an improvement of the convergence rate towards
  the optimal order.
This leads to the question on the limiting behavior. 
For the present example, $8M$ sweeps proved sufficient to approximate
  the latter.
Figure~\ref{fig:tg_sdc_m-v} shows the results of a corresponding study 
  with ${M = 1, \dots, 4}$ subintervals.
In contrast to Fig.~\ref{fig:tg_pdp_sdc_3-v} the limit cases attain 
  nearly constant slopes over a wide range of $\Delta t$.
With ${M\!=\!1}$ and $2$ the method achieves the optimal convergence
  rate of $2M$, whereas ${M\!=\!3}$ exhibits a slight order reduction
  of about $0.6$.
Comparing the latter to ${K\!=\!9}$ in Fig.~\ref{fig:tg_ppp_sdc_3-v} 
  reveals, however, that the error constant was reduced by nearly 
  two orders of magnitude.
With ${M\!=\!4}$ the method is still affected by the instability 
  of the explicit part for the two largest $\Delta t$.
After crossing the stability threshold it jumps almost instantly
  to the spatial error so that no asymptotic slope could be 
  determined.
In addition to the standard configuration, i.e. rotational velocity
  correction with divergence/mass-flux stabilization and final projection,
  Fig.~\ref{fig:tg_sdc_m-v} shows the results for ${M\!=\!3}$
  obtained using the standard velocity correction (${\chi\!=\!1}$)
  with or without stabilization, i.e. 
  ${\tau_\star\!=\!0}$ or
  ${\tau_\star\!=\!1}$, respectively,
  and no final projection.
They demonstrate the immense impact of the temporal splitting 
  error and violations of continuity that are caused by the 
  projection method used in the predictor and the corrector.
These errors prevent the convergence of the SDC method to the 
  underlying collocation method.
However, reducing the time step also diminishes the splitting and continuity 
  errors and may lead an apparent superconvergence as observed 
  in the case ${\chi\!=\!1}$ with ${\tau_\star\!=\!0}$.
Even the results of rotational method with stabilization and final projection
  may still differ from the Gauss collocation method.
Nevertheless, they meet the expected characteristics except for a mild 
  order reduction.
This is a substantial improvement over the SDPC method, which suffers a serious
  degradation of convergence for ${M>1}$, see
  \cite[Fig.~6.5]{Minion2018a_TI}.%

%

\subsection{Traveling 3D vortex with variable viscosity}


\subsubsection{Test cases}

The suitability for Navier-Stokes problems with variable viscosity
  is examined on the basis of a manufactured solution proposed by
  \cite{Niemann2018a}.
The exact velocity and pressure are given by
  \begin{alignat}{2}
    \label{eq:var-visc-test:vx}
    & v_x^{\mathrm{ex}}
    &&= \big[ \sin( 2\pi(x + t)) + \cos( 2\pi(y + t)) \big] \sin( 2\pi(z + t)) 
    \,,
    \\
    \label{eq:var-visc-test:vy}
    & v_y^{\mathrm{ex}}
    &&= \big[ \cos( 2\pi(x + t)) + \sin( 2\pi(y + t)) \big] \sin( 2\pi(z + t))
    \,,   
    \\
    \label{eq:var-visc-test:vz}
    & v_z^{\mathrm{ex}} 
    &&= \big[ \cos( 2\pi(x + t)) + \cos( 2\pi(y + t)) \big] \cos( 2\pi(z + t))
    \,,
    \\
    \label{eq:var-visc-test:p}
    & p^{\mathrm{ex}}   
    &&= \sin( 2\pi(x + t)) \sin( 2\pi(y + t)) \sin( 2\pi(z + t))
    \,.
  \end{alignat}
Equations (\ref{eq:var-visc-test:vx}\,--\,\ref{eq:var-visc-test:vz})
  define a periodic vortex array with wave length ${l=1}$ and
  velocity magnitude ${v^{\mathrm{ex}}_{\max} = 2}$, 
  traveling with a phase velocity of 1 in each direction, separately.
The exact solution is supplemented with a spatially and temporally 
  varying viscosity of the form
  ${\nu = \nu_0 + \nu_t(\V x, t, \V v)}$.
Three different scenarios are considered for the fluctuation $\nu_t$:
  \begin{alignat}{2}
    &\nu_t^{(1)}(\V x) 
    &&= \nu_1 \sin^2 (2\pi x) \sin^2 (2\pi y) \sin^2 (2\pi z)
    \,,
    \\[0.25\baselineskip]
    &\nu_t^{(2)}(\V x, t) 
    &&= \nu_1 \sin^2 (2\pi (x - t)) \sin^2 (2\pi (y - t)) \sin^2 (2\pi (z - t))
    \,,
    \\[0.25\baselineskip]
    &\nu_t^{(3)}(\V v) 
    &&= \nu_1 \frac{\V v^2}{|\V v^{\mathrm{ex}}|_{\max}^2}
    \, .
  \end{alignat}
Note that the expressions are normalized such that ${\max|\nu_t| = \nu_1}$ 
  provided that ${\V v = \V v^{\mathrm{ex}}}$.
The spatially varying fluctuation $\nu_t^{(1)}$ was already given in
  \cite{Niemann2018a}.
Complementing it with a unit phase velocity which is opposed to that of 
  $\V v^{\mathrm{ex}}$ leads to $\nu_t^{(2)}$\!\!.
Finally, $\nu_t^{(3)}$ depends on the approximate velocity and, thus,
  renders the viscous term genuinely nonlinear.
Based on the above specifications, the RHS of the Navier-Stokes problem
  is computed as
  \begin{equation}
    \label{eq:var-visc-test:rhs}
    \V f(\V x, t)
         = \d_t \V v^{\mathrm{ex}}
         + \nabla \cdot \V v^{\mathrm{ex}} \V v^{\mathrm{ex}}
         - \nabla \cdot \Big[ \nu(\V x, t, \V v^{\mathrm{ex}})
                              \Big( \nabla \V v^{\mathrm{ex}}
                                  + \transpose{(\nabla \V v^{\mathrm{ex}})}
                              \Big)
                         \Big]
         + \nabla p^{\mathrm{ex}}
    \, .
  \end{equation}
Omitting the convection term in \eqref{eq:var-visc-test:rhs} yields
  the RHS of the corresponding Stokes problem. 

In all studies reported below time-dependent Dirichlet conditions are
  imposed such that 
  ${\V v_{\mathrm b} = \V v^{\mathrm{ex}}}$ 
  on $\d\Omega$.


\subsubsection{Influence of the base time-integration scheme}

\begin{figure}
  \includegraphics[height=67mm]{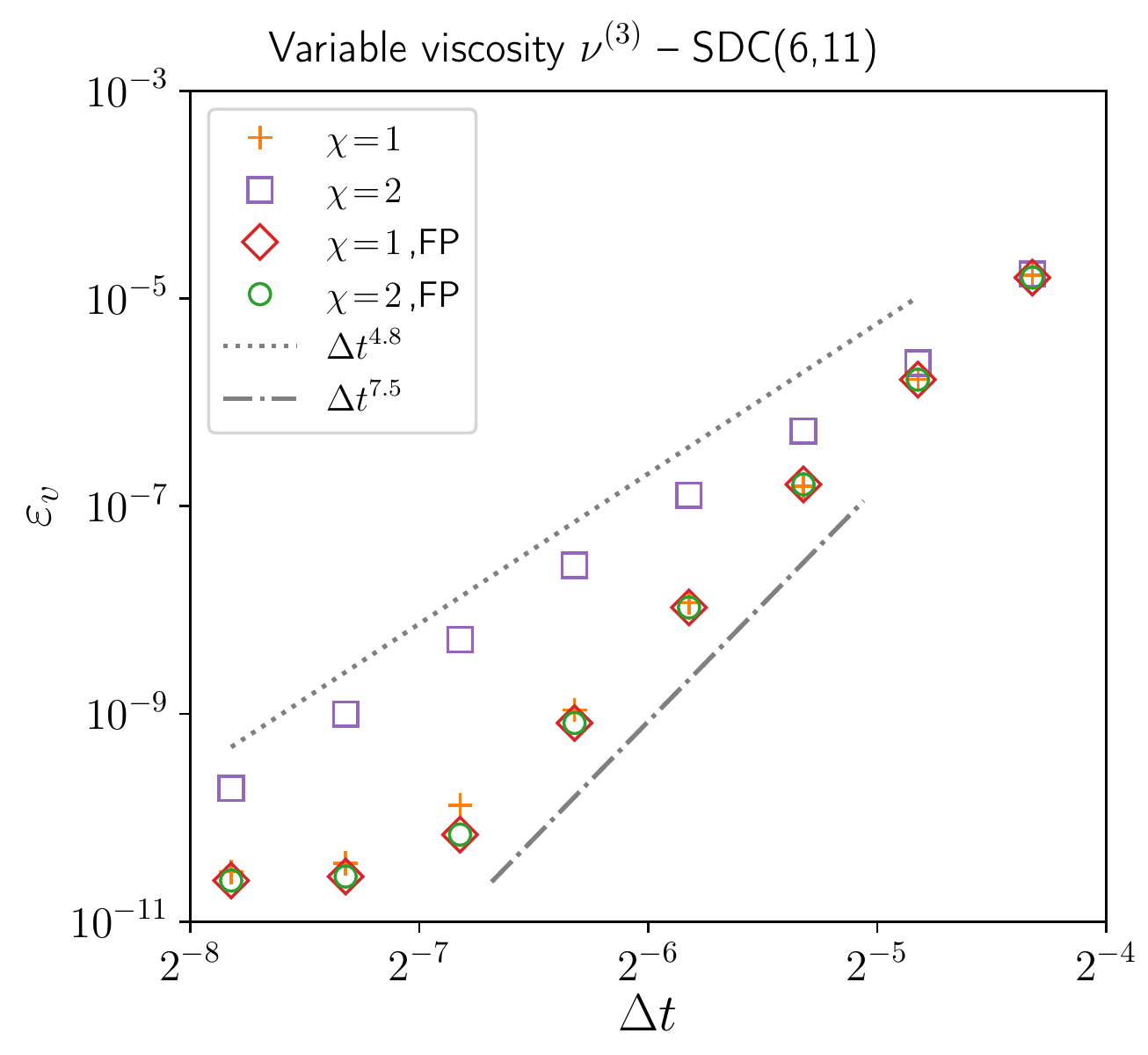}
\caption{%
  Velocity error of SDC with ${M=6}$ subintervals and ${K=11}$ correction
  sweeps for the traveling 3D vortex with solution-dependent viscosity.
  The choices
  ${\chi\!=\!1}$ 
  and
  ${\chi\!=\!2}$ 
  refer to the method based on standard and rotational
  velocity correction, respectively, and 
  FP indicates the application of the final projection step.
  \label{fig:vv_ddd_sdc_6-11_methods}
  }
\end{figure}

To investigate the role of the parameter $\chi$ and the final projection 
  in the predictor and corrector, a preliminary study was conducted
  for a solution-dependent viscosity $\nu^{(3)}$ with coefficients
  ${\nu_0 = \nu_1 = 10^{-2}}$.
This choice corresponds to Reynolds number of
  ${\mathit{Re} ={l v_{\max}}/{\nu_{\mathrm{ref}}} \approx 133}$,
  where
  ${\nu_{\mathrm{ref}} = \nu_0 + \frac{1}{2}\nu_1}$.
The numerical tests were computed in the domain 
  ${\Omega = [-\sfrac{1}{2},\sfrac{1}{2}]^3}$
  for
  ${t \le 0.25}$.  
Spatial discretization is based on a uniform mesh comprising 
  $2^3$ cubic elements of degree ${P=16}$. 
For time integration the SDC method was applied with ${M=6}$ and ${K=11}$
  correction sweeps.
Figure~\ref{fig:vv_ddd_sdc_6-11_methods}
  shows the resulting velocity error for different predictor/corrector
  variants.
As in case of the Taylor-Green vortex, the variants with FP
  achieve the best results and virtually coincide regardless of the choice
  for $\chi$. 
They reach a convergence rate of approximately 7.5, which corresponds 
  to a reduction of 4.5 or 37.5 percent of the expected order of 12.
The standard scheme (${\chi=1}$) with no FP attains almost the 
  same accuracy, whereas the rotational scheme (${\chi=2}$) converges
  at a rate of only about 4.8.
This contradicts the results obtained with the Taylor-Green vortex, 
  for which the rotational scheme surpassed the standard one. 
These observations indicate that the final projection eliminates a 
  substantial part of the splitting error, while the choice of $\chi$ 
  is of minor importance.
Based on these observations, all of the following studies use 
  the rotational scheme with FP.

\subsubsection{Temporal convergence}

\begin{figure}
\subfloat[Constant viscosity.]
  {\includegraphics[height=67mm]{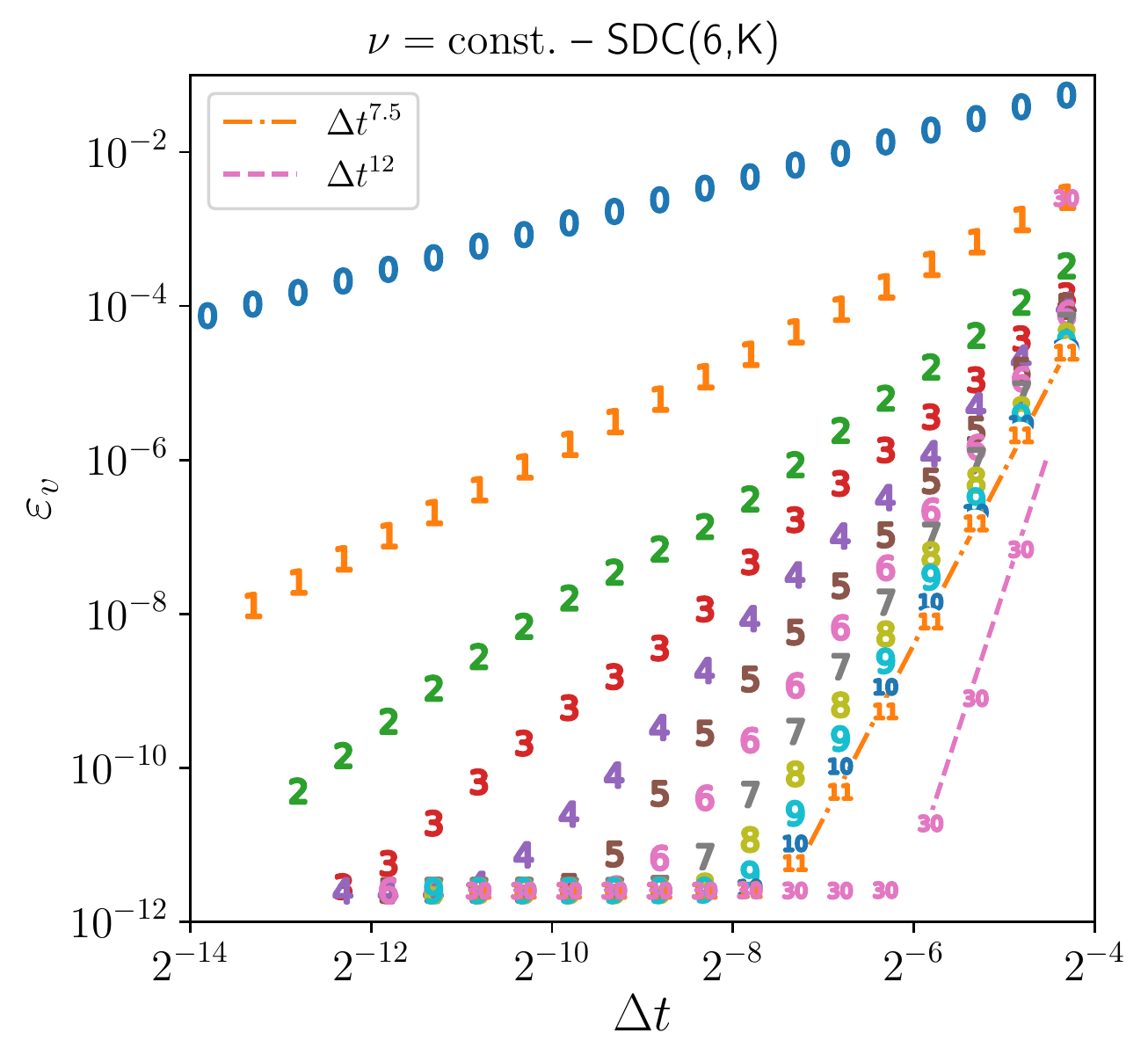}
   \label{fig:vv0_ddd_sdc_6-v}}
\subfloat[Spatially varying viscosity.]
  {\includegraphics[height=67mm]{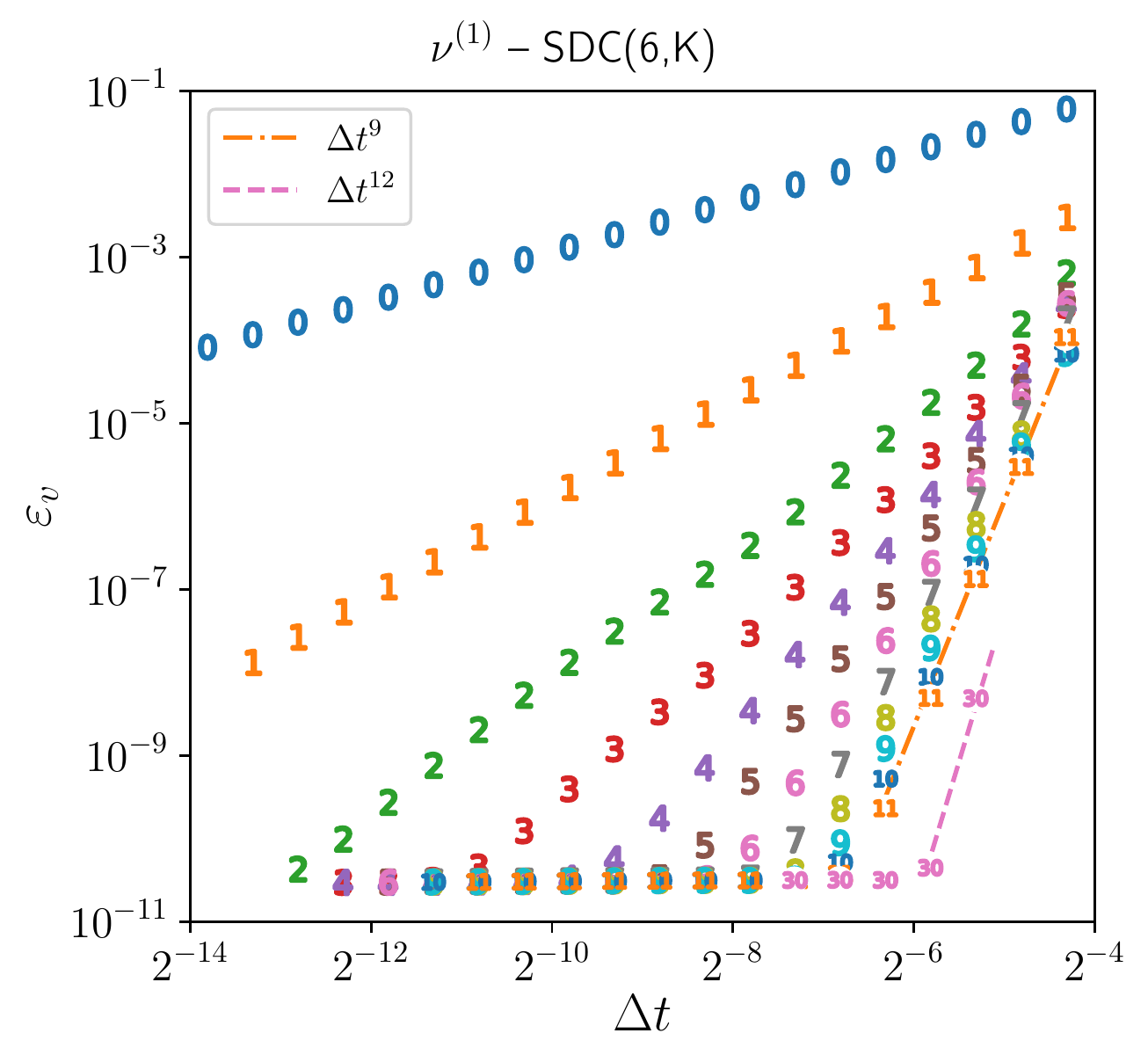}
   \label{fig:vv1_ddd_sdc_6-v}}
\\
\subfloat[Spatiotemporally varying viscosity.]
  {\includegraphics[height=67mm]{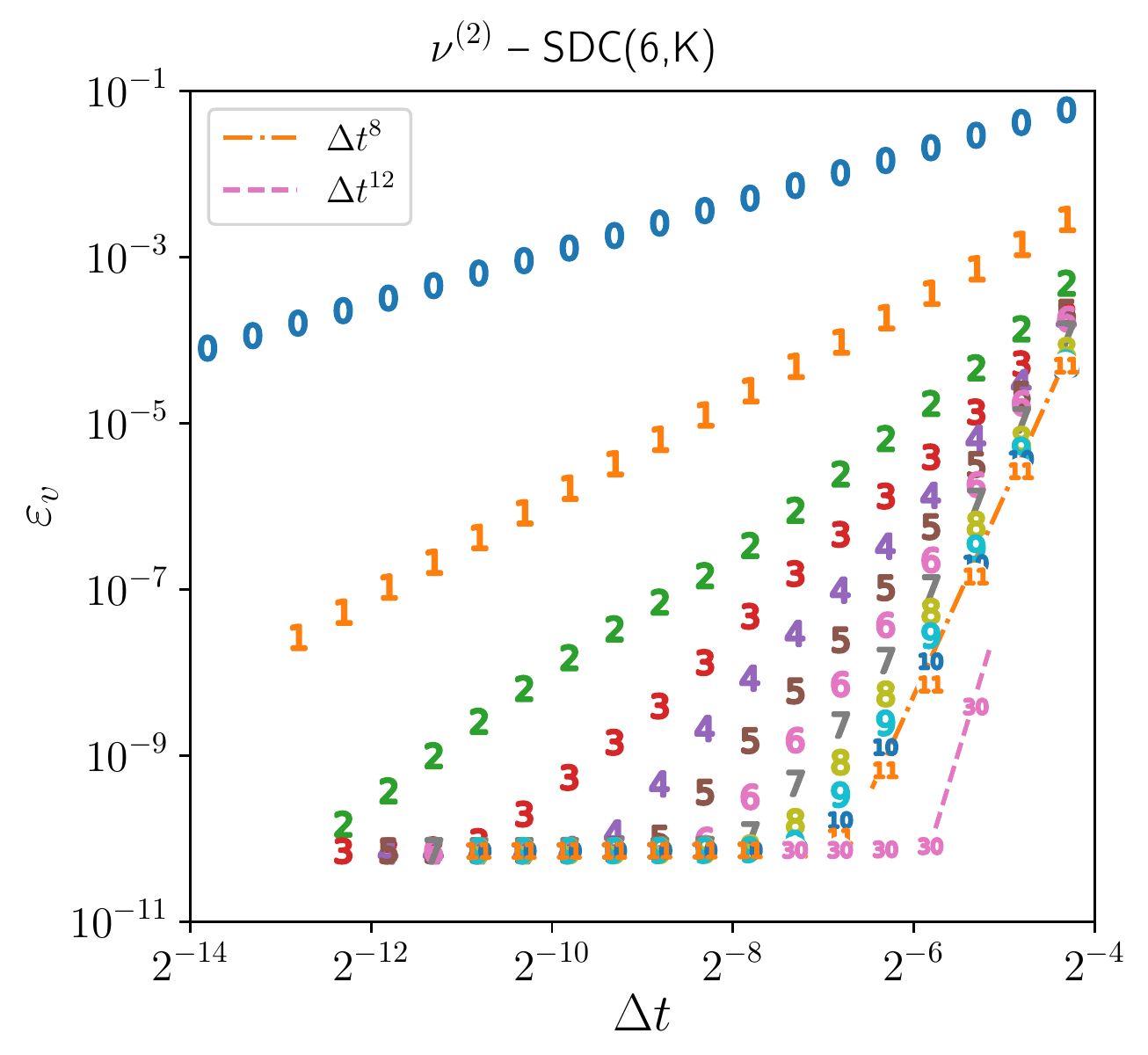}
   \label{fig:vv2_ddd_sdc_6-v}}
\subfloat[Solution-dependent viscosity.]
  {\includegraphics[height=67mm]{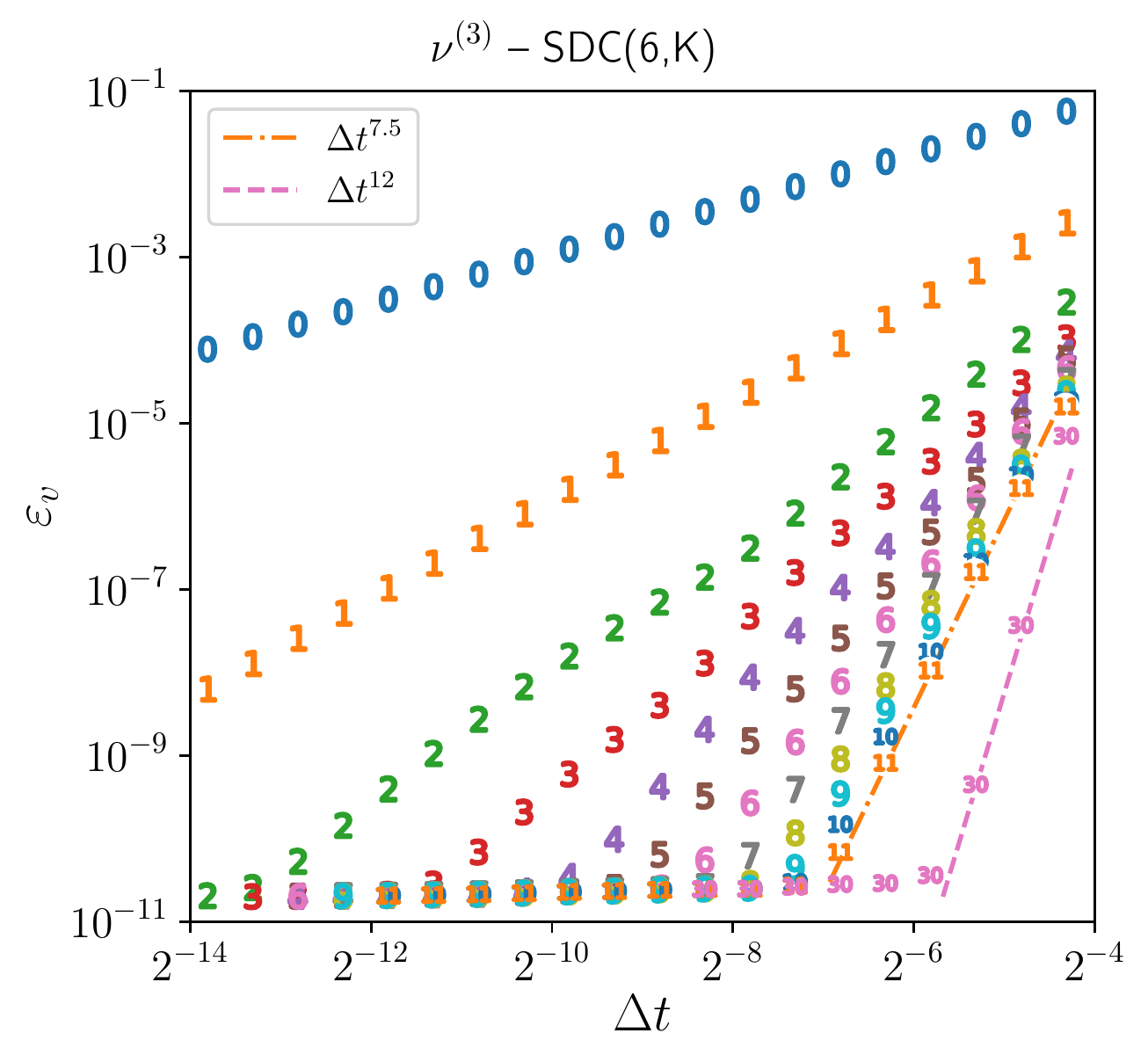}
   \label{fig:vv3_ddd_sdc_6-v}}
\caption{%
  Velocity error for different viscosity scenarios using
  SDC with ${M=6}$ subintervals.
  Labels indicate the number of correction sweeps and
  lines the approximate slope.
  \label{ddd_sdc_6_v}
  }
\end{figure}

Three different scenarios were chosen for assessing the robustness of 
  the SDC method against variable viscosity:
  1) spatially varying,          ${\nu(\V x)    = \nu^{(1)}}$, 
  2) spatiotemporally varying,   ${\nu(\V x, t) = \nu^{(2)}}$, and
  3) solution-dependent,         ${\nu(\V v)    = \nu^{(3)}}$,
  with coefficients 
  ${\nu_0 = \nu_1 = 10^{-2}}$.
Additionally, the case with constant ${\nu = 0.015}$ is considered for reference.
The computational domain, spatial discretization and final time are identical
  to the previous study.
Figure~\ref{ddd_sdc_6_v} shows the velocity error obtained with ${M=6}$
  subintervals and a different number of correction sweeps, ranging from
  ${K=0}$ up to $30$.
All investigated scenarios exhibit a similar behavior and achieve convergence
  rates comparable to the reference case.
Using 30 corrections yields a rate of 12, which equals the theoretical order 
  of the underlying collocation method.
Runs with a lower number of sweeps suffer an order reduction.
The extent of this reduction is similar for all scenarios, which indicates 
  that the presence of a variable viscosity is not the primary cause.
As in the previous test case it is more likely to be caused by the boundary
  treatment and the splitting error of the underlying projection method.
It is further noted that 
  Fig.~\ref{fig:vv1_ddd_sdc_6-v} and \ref{fig:vv2_ddd_sdc_6-v}
  lack the errors for the two largest time steps with ${K=30}$.
This is because $\Delta t$ exceeds the long term stability threshold,
  which is considered in more detail below.

\begin{figure}
\subfloat[Pressure error.]
  {\includegraphics[height=67mm]{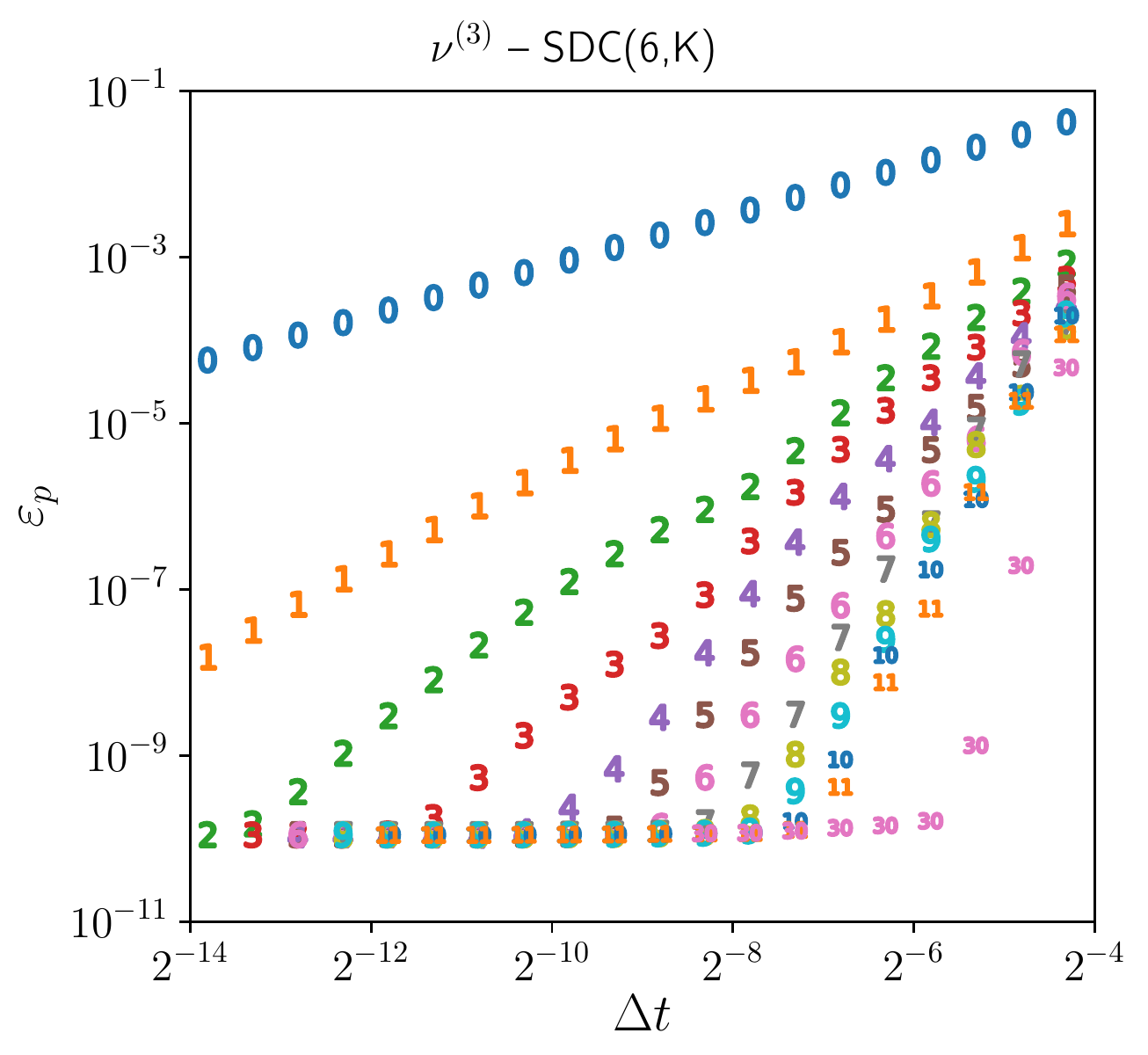}
   \label{fig:vv3_ddd_sdc_6-p}}
\subfloat[Divergence.]
  {\includegraphics[height=67mm]{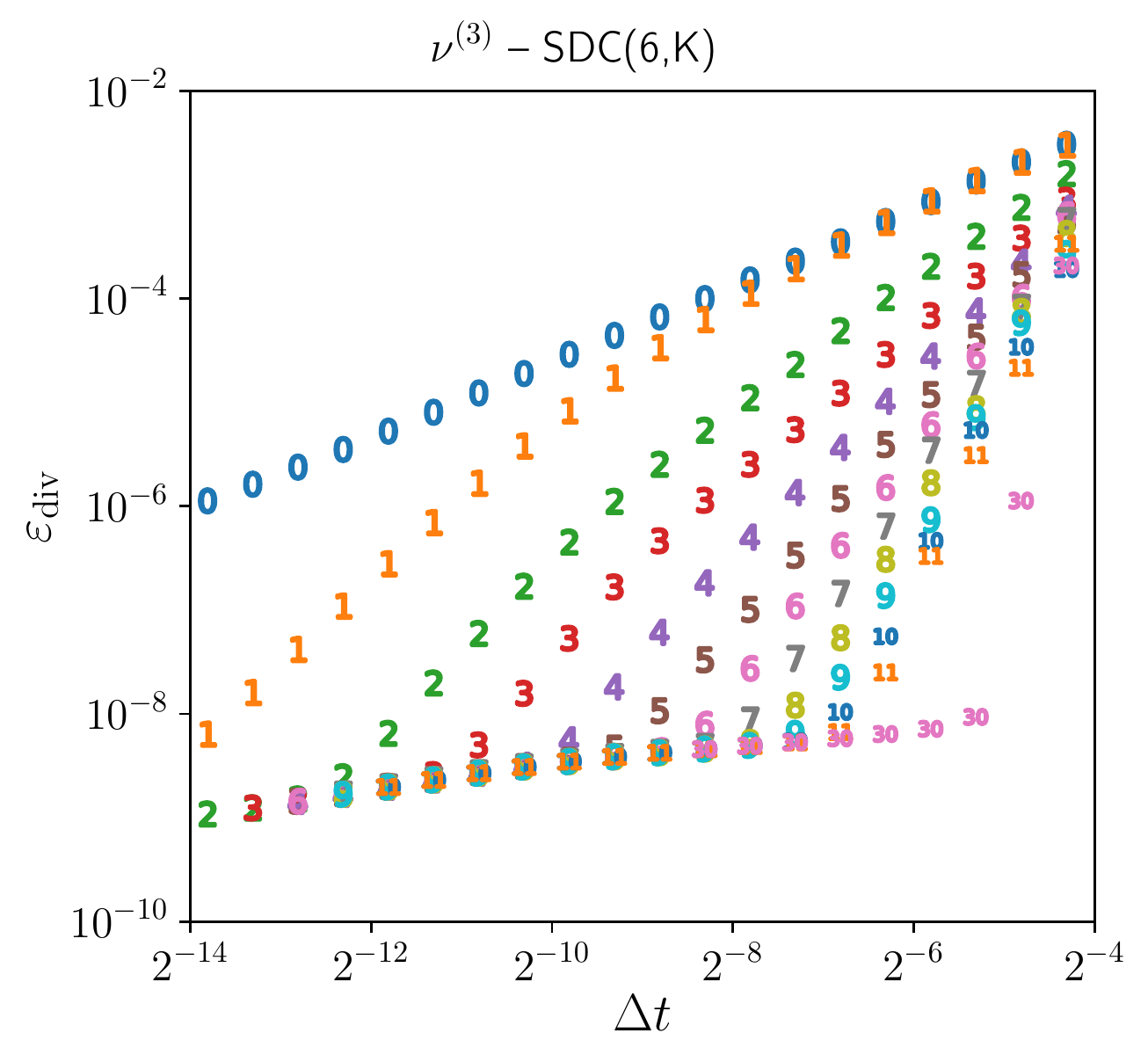}
   \label{fig:vv3_ddd_sdc_6-div}}
\caption{%
  Pressure and divergence errors obtained
  for the case of solution-dependent viscosity.
  \label{ddd_sdc_6_p-div}
  }
\end{figure}

As discussed in Sec.~\ref{sec:sdc:flow:corrector}, the SDC method does 
  not provide the discrete pressure.
It can be computed, however, by solving the discrete version of the 
  consistent pressure equation \eqref{eq:p:consistent:pde}.
This yields the pressure with an accuracy comparable to that of velocity,
  see, e.g., Fig.~\ref{fig:vv3_ddd_sdc_6-p} for the case of solution-dependent
  viscosity.
Figure~\ref{ddd_sdc_6_p-div} depicts the corresponding divergence error.
Except for larger time steps with ${K=1}$ it shows roughly the same behavior
  as the velocity and pressure errors.
Additional studies revealed that an even stronger divergence penalization
  fails to reduce $\varepsilon_{\mathrm{div}}$ significantly.
This observation is somewhat surprising.
It can be explained with the splitting
  error in the final projection step, which incurs a violation of the 
  tangential velocity boundary conditions and causes a growth of divergence 
  near the edges of the computational domain.

\subsubsection{Convergence towards the collocation solution}

\begin{figure}
  \includegraphics[height=67mm]{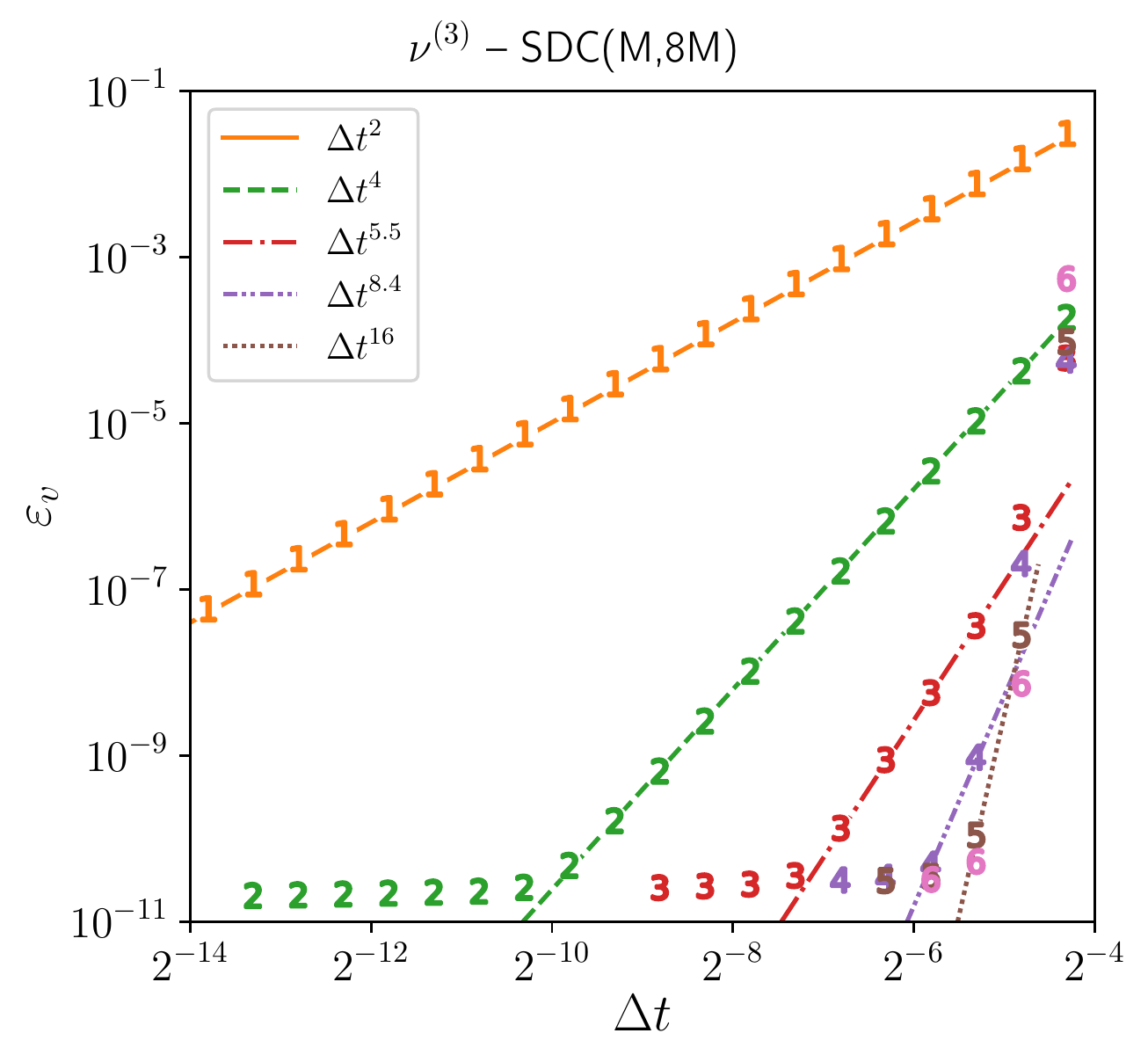}
\caption{%
  Velocity error of the converged SDC method for the case of
  solution-dependent viscosity. 
  Labels indicate the number $M$ of subintervals.
  \label{fig:ddd_sdc_6_m-v}
  }
\end{figure}

To assess the limiting behavior for variable viscosity, the SDC method was
  applied with ${M = 1, \dots, 6}$ subintervals and ${K\!=\!8M}$ correction
  sweeps.
Figure~\ref{fig:ddd_sdc_6_m-v} shows the most challenging case with a
  solution-dependent viscosity.
Similarly to the Taylor-Green example, regular convergence with the expected 
  rate of $2M$ is observed for ${M\!=\!1}$ and $2$, whereas ${M\!=\!3}$ exhibits
  a slight reduction of half an order.
For time steps ${\Delta t > 2^{-5}}$ the method is affected by the instability
  of the explicit part with $M \ge 3$.
Additionally, the splitting error may still be significant here. 
The cumulative effect of both factors could explain the superconvergence
  observed for $M \ge 4$.

\subsubsection{Stability}

The stability of the SDC method was investigated for the limiting cases of
  convection-dominated flow and Stokes flow.
For this study the domain 
  ${\Omega = [-\sfrac{1}{2},\sfrac{1}{2}]^3}$
  was discretized in three ways:
  1) $6^3$ elements of degree ${P\!=\!6}$,
  2) $3^3$ elements of degree 11 and
  3) $2^3$ elements of degree 16.
The numerical tests were run until reaching the final time ${T = 10}$, 
  which corresponds to 10 convective units in terms of the phase velocity.
A test was considered unstable when 
  exceeding a velocity magnitude of ${v^{\mathrm{ex}}_{\max}/2}$ or
  detecting a NaN in the numerical solution.
Based on this criterion the time step was adapted via bisection until
  a reaching sufficiently accurate approximation of the critical 
  time step $\Delta t^\ast$.

\begin{figure}
  \includegraphics[height=67mm]{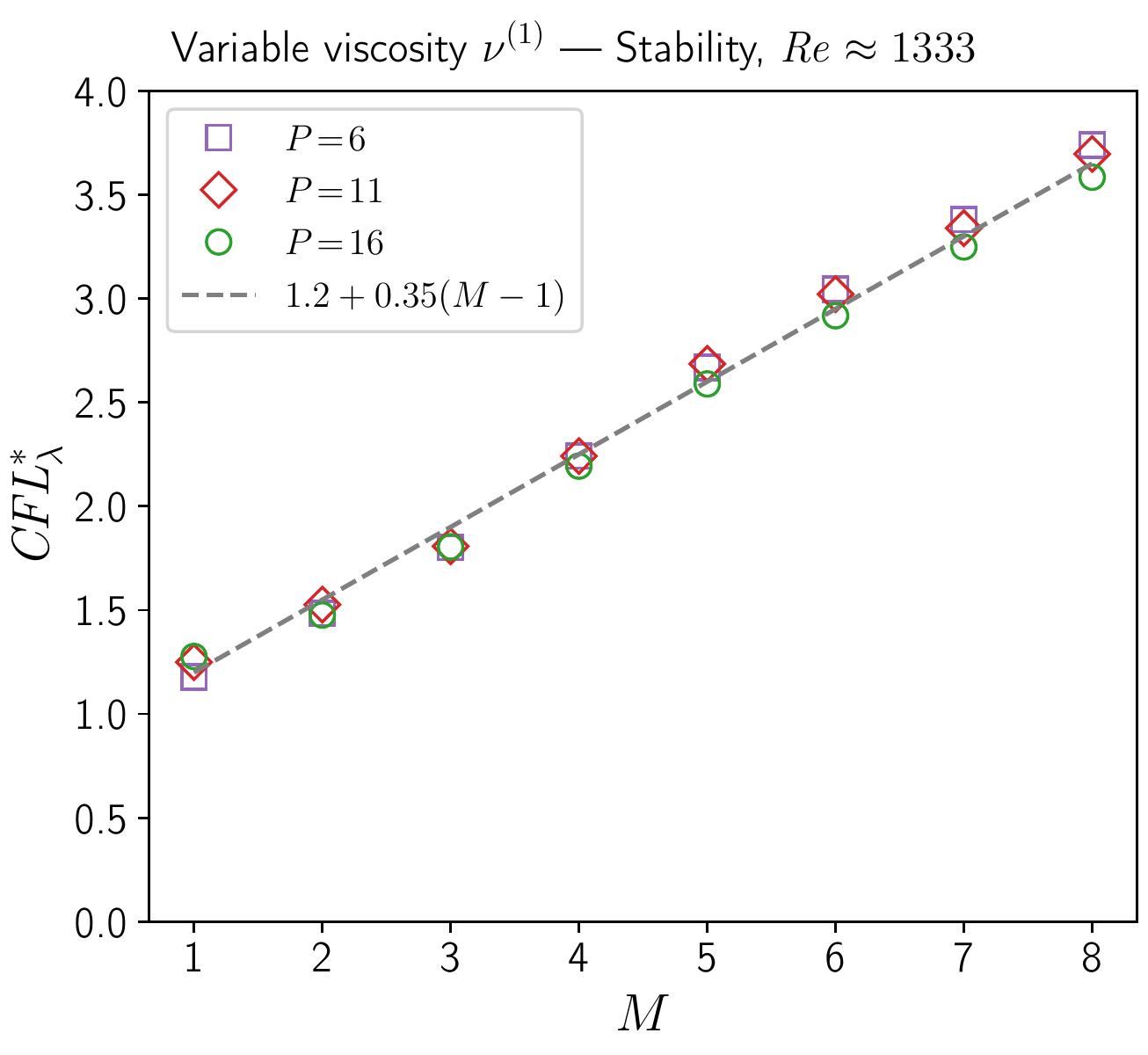}
  \caption{%
    Critical CFL number for a varying number of subintervals $M$ 
    and ${K=2M}$ correction sweeps.
    \label{fig:vv1_ddd_cfl_lambda}
    }
\end{figure}

In the convection-dominated case the viscosity $\nu^{(1)}$ is adopted
  with ${\nu_0 = \nu_1 = 10^{-3}}$.
This corresponds to a Reynolds number of approximately 1333,
  based on $v^{\mathrm{ex}}_{\max}$ and
  ${\nu_{\mathrm{ref}} = \nu_0 + \frac{1}{2} \nu_1}$.
The resulting convective stability threshold is converted into dimensionless form 
  by introducing the critical Courant-Friedrichs-Lewy (CFL) number
  \begin{equation}
    \CFLc = \frac{\Delta t^\ast v_{\max}}{\delta}
    \,,
  \end{equation}
  where $\delta$ is a length scale characterizing the mesh spacing.
For high-order element methods different length scales have been proposed,
  in particular,
  ${\delta_1 = h / (P+1)}$ and
  ${\delta_2 = h / P^2}$, 
  see e.g. \cite{Gassner2011a_SE,Karniadakis2005a_SE}.
The corresponding CFL numbers are denoted as $\CFLc_1$ and $\CFLc_2$, respectively.
Alternatively, the length scale can be defined as the inverse maximum modulus 
  of the eigenvalues of the one-dimensional element convection operator, i.e.,
  ${\delta_{\lambda} = |\lambda|_{\max}^{-1}}$.
In the present case, the eigenvalues result from the GLL collocation differentiation
  operator of degree $P$ combined with one-sided Dirichlet conditions.
For details see \citet[Sec. 7.3.3]{Canuto2011a_SE}.
The resulting CFL number is denoted as $\CFLcLmb$.
Tabular~\ref{tab:cfl} compiles the critical CFL numbers determined in this study.
Note that the first three lines correspond to the case with only one subinterval.
Starting with the split semi-implicit Euler method (${K\!=\!0}$), the admissible
  time step grows by factor of 3 with one correction and even a factor 
  of $5.3$ with two. 
Elevating the number of subintervals and, proportionally, the number correction 
  sweeps yields a further increase of the stability threshold.
A comparison of the stability limits for equal $M$ further reveals 
  a strong dependence of $\CFLc_1$ and $\CFLc_2$ on $P$,
  whereas $\CFLcLmb$ is virtually independent of the polynomial degree.
Figure~\ref{fig:vv1_ddd_cfl_lambda} confirms this observation and, moreover,
  illustrates that the stability threshold grows linearly when increasing the 
  number of subintervals.
This behavior was expected, since the predictor and the corrector perform 
  substeps with the size scaling as ${\Delta t/M}$.
Given the influence of
  the flow configuration,
  the Reynolds number,
  the stability and dissipativity of spatial discretization
  and 
  the termination criterium,
  these results cannot be compared directly to other studies.
Despite this limitation, the $\CFLc_2$ observed for 
  the second-order SDC method (${M\!=\!1}$, ${K\!=\!2}$)
  reside in the same range 
  as those reported by \citet{Fehn2017a_TI},
  who used a similar space discretization combined with a 
  semi-implicit projection method of order 2.

It should be noted that the SDC method based on IMEX Euler is potentially 
  unstable in the purely convective case.
For example, according to \cite[Fig.~4.4]{Minion2003b_TI}
  the stability region includes no part of the imaginary axis
  for $K=M$ with $M=5$ and $6$.
Moreover, the convection term of the flow problem may give rise to
  nonlinear instability.
The arising instabilities often grow with further correction sweeps
  so that no convergence to the corresponding Gauß method can be
  achieved.
Both effects, marginal linear stability and nonlinearity, prevent the 
  deduction of a reliable stability criterion for convection, 
  which leaves an inconvenience for the application to high Reynolds number
  flows.

\begin{table}
  \caption{Critical CFL numbers for ${\mathit{Re} \approx 1333}$ 
    and SDC with a different number of subintervals $M$ and corrections sweeps $K$.
    \label{tab:cfl}}
  \begin{tabular}{cccccccccccccc} \toprule
       &      && \multicolumn{3}{c}{$P=6$}          && \multicolumn{3}{c}{$P=11$}         && \multicolumn{3}{c}{$P=16$} \\ \cmidrule{4-6}\cmidrule{8-10}\cmidrule{12-14}
   $M$ & $ K$ && $\CFLc_1$ & $\CFLc_2$ & $\CFLcLmb$ && $\CFLc_1$ & $\CFLc_2$ & $\CFLcLmb$ && $\CFLc_1$ & $\CFLc_2$ & $\CFLcLmb$ \\ \midrule
   $1$ & $ 0$ && $.174$    & $.895$    & $.222$     && $.142$    &  $1.43$   &  $.277$    &&  $.118$   &  $1.78$   &   $.317$   \\
   $1$ & $ 1$ && $.522$    & $2.68$    & $.667$     && $.368$    &  $3.72$   &  $.720$    &&  $.275$   &  $4.14$   &   $.737$   \\
   $1$ & $ 2$ && $.922$    & $4.74$    & $1.18$     && $.639$    &  $6.45$   &  $1.25$    &&  $.476$   &  $7.17$   &   $1.28$   \\
   $2$ & $ 4$ && $1.16$    & $5.97$    & $1.48$     && $.781$    &  $7.87$   &  $1.53$    &&  $.551$   &  $8.29$   &   $1.48$   \\
   $3$ & $ 6$ && $1.41$    & $7.25$    & $1.80$     && $.925$    &  $9.32$   &  $1.81$    &&  $.673$   &  $10.1$   &   $1.80$   \\
   $4$ & $ 8$ && $1.75$    & $9.02$    & $2.24$     && $1.14$    &  $11.6$   &  $2.24$    &&  $.818$   &  $12.3$   &   $2.19$   \\
   $5$ & $10$ && $2.09$    & $10.7$    & $2.67$     && $1.37$    &  $13.8$   &  $2.68$    &&  $.966$   &  $14.5$   &   $2.59$   \\
   $6$ & $12$ && $2.38$    & $12.3$    & $3.04$     && $1.54$    &  $15.6$   &  $3.02$    &&  $1.09$   &  $16.4$   &   $2.91$   \\
   $7$ & $14$ && $2.64$    & $13.6$    & $3.38$     && $1.71$    &  $17.2$   &  $3.34$    &&  $1.21$   &  $18.2$   &   $3.25$   \\
   $8$ & $16$ && $2.92$    & $15.0$    & $3.74$     && $1.89$    &  $19.1$   &  $3.70$    &&  $1.34$   &  $20.1$   &   $3.58$   \\ \bottomrule
  \end{tabular}
\end{table}

A similar study was conducted for the Stokes flow with
  solution-dependent viscosity $\nu^{(3)}$ featuring
  ${\nu_0 = \nu_1 = 1}$.
Notwithstanding the semi-implicit treatment of the viscous term
  all test runs remained stable up to the maximal time step size
  of ${\Delta t = 5}$, which is about $10^4$ times the 
  diffusive fluctuation time scale 
  ${\tau_{\mathrm d} = \sfrac{h^2}{\nu_1 P^{4}}}$.
This allows the conclusion that the semi-implicit approach does not affect 
  the stability of the SDC method for viscosity fluctuations up to at least 
  50 percent.


\subsubsection{Spatial convergence}

\begin{figure}
  \includegraphics[height=67mm]{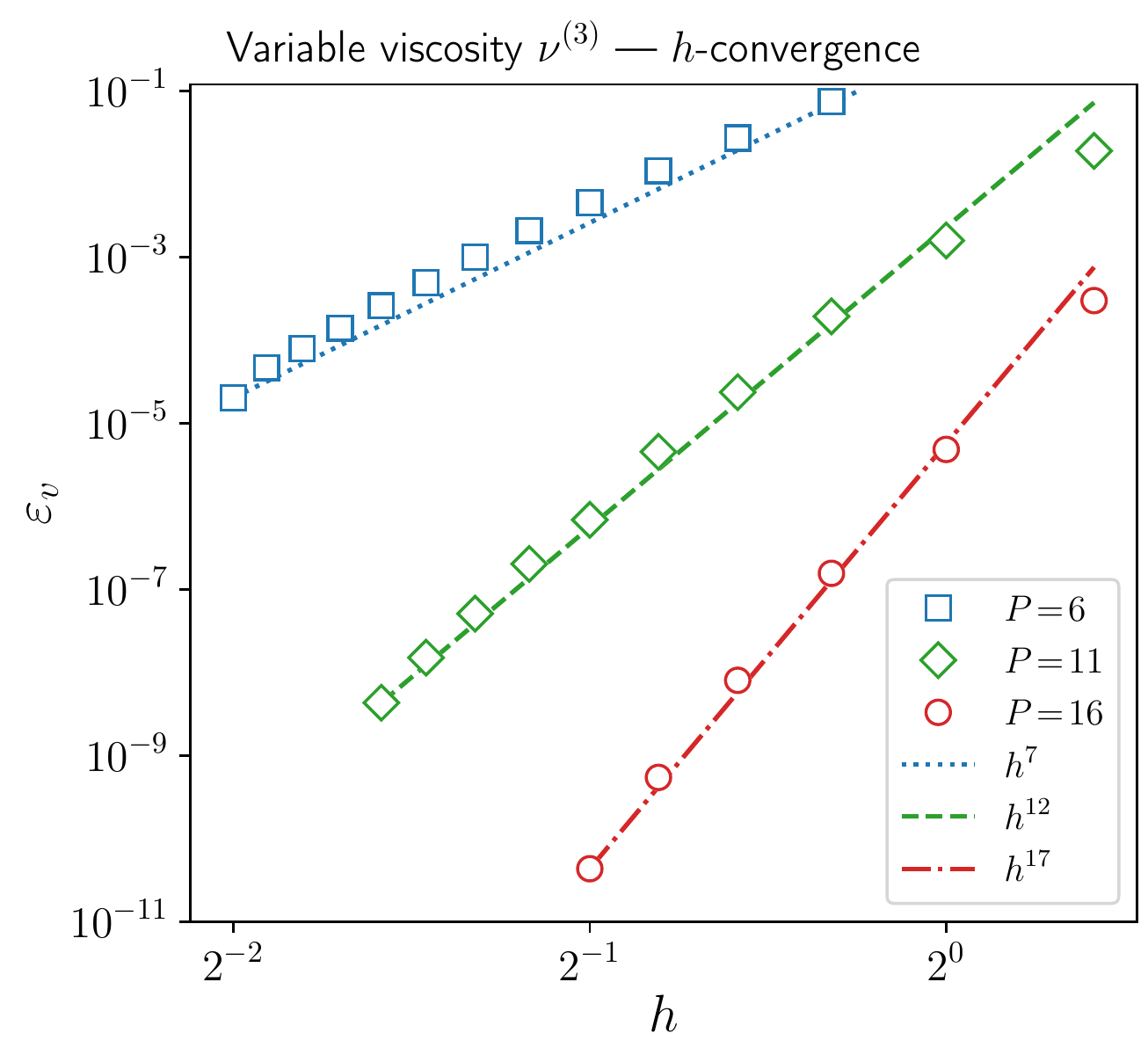}
  \caption{%
    Spatial convergence with solution-dependent variable viscosity. 
    \label{fig:vv3_h-convergence}
    }
\end{figure}

The final study serves for examining the spatial convergence with variable 
  viscosity.
It is based on the on the highly-nonlinear, solution-dependent viscosity 
  $\nu^{(3)}$ with coefficients ${\nu_0 = \nu_1 = 10^{-2}}$. 
The numerical tests were computed in the domain ${\Omega=[-2,2]^3}$
  for polynomial degrees up to ${P=16}$.
Starting from one, the number of elements per direction was gradually increased
  to at least 8.
Time integration was performed using the SDC method 
  with ${M=6}$ and ${K=11}$ until reaching ${T=0.25}$.
The time step size was confined to ${\Delta t = 2^{-8}}$ 
  such that the temporal discretization error is negligible.
Figure~\ref{fig:vv3_h-convergence} shows the velocity error for 
  polynomial degrees ${P=6}$, $11$ and $16$.
In all three cases the method converges approximately with $h^{P+1}$.
This result is surprising, since the viscous terms are integrated
  with just ${P+1}$ GLL points, for which only order $P$ is expected
  \cite{Maday1990a_SE}.
Possibly, the higher convergence rate is promoted by the construction of 
  the test case.
Clarifying this issue requires an in-depth investigation, which is 
  beyond the scope of the present work.

\section{Conclusions}
\label{sec:conclusions}

This paper presents a high-order, semi-implicit time-integration strategy  
  for incompressible Navier-Stokes problems with variable viscosity based 
  on the spectral deferred correction (SDC) method.
Combining SDC with the discontinuous Galerkin spectral-element 
  method for spatial discretization yields a powerful
  approach for targeting arbitrary order in space and time.
The key ingredients of the method, the predictor and the corrector, are derived
  from a first-order velocity-correction method, which is augmented by 
  an additional projection step to remove divergence errors caused
  by variable viscosity.
In contrast to the SDPC method of \citet{Minion2018a_TI}, the pressure 
  occurs only as an auxiliary variable in the substeps and needs not
  to be stored.
Furthermore, mixed-order (${P/P\!-\!1}$) polynomial approximations are used for 
  inf-sup stability and combined with divergence/mass-flux stabilization 
  for pressure robustness \cite{Akbas2018a_SE}.

The performance of the SDC method was assessed at the example of a Taylor-Green
  (TG) vortex and a manufactured 3D vortex array, both traveling with a prescribed
  phase velocity.
The latter involves a variable viscosity $\nu$ that can be chosen to depend on
  space,
  space and time, or on
  the velocity $\V v$.
For the TG vortex with constant viscosity and periodic boundaries, each correction
  sweep lifted the order by 1 as expected.
Imposing time-dependent Dirichlet boundary conditions, however, led to a more
  irregular convergence behavior accompanied by an order reduction.
For a fixed number of corrections, convergence starts at a lower rate and at a
  higher error level than in the periodic case.
With decreasing size of the time step $\Delta t$ it accelerates and can even 
  reach a rate higher than the theoretical order.
A closer inspection revealed that this behavior is most likely a manifestation
  of the splitting error of the time integration scheme adopted for the corrector.
Using the rotational velocity correction augmented by divergence/mass-flux
  stabilization and a final projection step allowed to minimize this error
  and lowered the resulting order reduction to approximately $0.2$ per sweep.
The investigation showed that the splitting error also affects the limit case
  which arises after a sufficient number of corrections.
While the optimized corrector based on the rotational velocity-correction
  exhibited a nearly ideal convergence behavior with only a slight order 
  reduction, other variants suffered a severe degradation.
Moreover, these observations offer an explanation for the results of 
  \citet{Minion2018a_TI} who achieved only about half the expected 
  order in the limit case.
Thus, it is not surprising that the proposed SDC method clearly outperforms 
  SDPC except for the periodic case, for which the splitting error disappears.
Finally, the new method proved robust and almost equally efficient with
  variable viscosity, even in the nonlinear case, where the latter depends
  on the solution itself.
This case also served as a test bed for exploring the temporal stability.
Considering a convection-dominated regime, the critical CFL number grows 
  linearly with the number of SDC subintervals and virtually coincides for 
  different polynomial degrees, when choosing the length scale based on 
  eigenvalues of the DG differentiation operator.
In the Stokes case with variable viscosity, the method remains stable for 
  steps up to at lest 4 orders of magnitude above the viscous time scale 
  $\sfrac{h^2}{P^4 \nu}$, where $h$ is the element size and $P$ the 
  polynomial degree of the discrete velocity.

The present study suggests that a further improvement is possible by 
  eliminating the residual splitting error.
This could be achieved by performing a coupled iteration within each Euler
  step or using a higher order method for constructing the corrector.
Of course, the ultimate goal is to render the SDC method competitive 
  to common approaches such as multistep and Runge-Kutta methods.
This may require further measures.
For example, a higher-order time-integration method can be harnessed 
  as the predictor to gain a more accurate initial approximation
  and conditional stability in the convective limit
  \cite{Layton2007a_TI}.
Using uniform instead of Gaussian quadrature points in time, such methods 
  may also accelerate the corrector \cite{Christlieb2015a_TI}.
Alternatively, the number of iterations could be reduced by adopting
  suitable preconditioners like the LU-decomposition proposed in
  \cite{Weiser2015a_TI}
  or
  space-time multilevel strategies such as MLSDC \cite{Speck2015a_TI}
  and their parallel descendants \cite{Minion2015a_TI,Speck2018a_TI}.
Utilizing these techniques for designing faster flow solvers is the 
  topic of ongoing research and will be addressed in future papers.
%

  
\section*{Acknowledgements}

Funding by German Research Foundation (DFG) in frame of the project 
  STI 57/8-1 is gratefully acknowledged.
The author would like to thank ZIH for providing computational resources
  and Karl Schoppmann for his assistence in devising and implementing the
  manufactured solution for the variable viscosity test case.



\end{document}